\documentclass{article}
\usepackage{amssymb,amsmath,bm,geometry,graphics,graphicx,url,color}

\def\no{\noindent}
\def\pmatrix{\left(\begin{array}}
\def\endpmatrix{\end{array}\right)}

\newtheorem{mydef}{Definition}[section]
\newtheorem{mytheo}[mydef]{Theorem}
\newtheorem{myrem}[mydef]{Remark}

\newtheorem{assum}[mydef]{Assumption}

\newtheorem{prop}[mydef]{Proposition}

\newcommand{\norm}[1]{\left\Vert#1\right\Vert}
\newcommand{\abs}[1]{\left\vert#1\right\vert}

\def\sinc{\mathrm{sinc}}

\title{A  long-term numerical energy-preserving analysis of symmetric and/or symplectic extended RKN integrators for efficiently solving highly oscillatory Hamiltonian systems}

\author{Bin Wang\,
\footnote{School of Mathematical Sciences, Qufu Normal University,
Qufu 273165, P.R. China; Mathematisches Institut, University of
T\"{u}bingen, Auf der Morgenstelle 10, 72076 T\"{u}bingen, Germany.
The research is supported in part by the Alexander von Humboldt
Foundation and by the Natural Science Foundation of Shandong
Province (Outstanding Youth Foundation) under Grant ZR2017JL003.
E-mail:~{\tt wang@na.uni-tuebingen.de} } \and Xinyuan
Wu\thanks{School of Mathematical Sciences, Qufu Normal University,
Qufu 273165, P.R. China; Department of Mathematics, Nanjing
University, Nanjing 210093, P.R. China. The research is supported in
part by the National Natural Science Foundation of China under Grant
11671200. E-mail:~{\tt xywu@nju.edu.cn}} }

\begin{document}
\maketitle

\begin{abstract} The primary objective of this paper is to present a long-term
numerical energy-preserving analysis of one-stage explicit symmetric
and/or symplectic extended Runge--Kutta--Nystr\"{o}m (ERKN)
integrators for highly oscillatory Hamiltonian systems. We study the
long-time numerical energy conservation not only for symmetric
integrators but also for symplectic integrators. In the analysis, we
neither assume symplecticity for symmetric methods, nor assume
symmetry for symplectic methods. It turns out that these both kinds
of ERKN integrators have a near conservation of the total and
oscillatory energy over a long term. To prove the result for
symmetric integrators, a relationship between symmetric ERKN
integrators and trigonometric integrators is established by using
Strang splitting and  based on this connection, the long-time
conservation is derived. For the long-term analysis of symplectic
ERKN integrators, the above approach does not work anymore and we
use the technology of modulated Fourier expansion developed in SIAM
J. Numer. Anal. 38 (2000) by Hairer and Lubich. By taking some novel
adaptations of this essential technology for non-symmetric methods,
we derive the modulated Fourier expansion for symplectic ERKN
integrators. Moreover, it is shown that the symplectic ERKN
integrators  have two almost-invariants and then the
 near energy conservation over a long term is obtained.
\medskip

\no{\bf Keywords:} Long-time  energy conservation\and Modulated
Fourier expansions\and Symmetric or symplectic methods\and Extended
RKN integrators\and Highly oscillatory Hamiltonian systems

\medskip
\no{\bf MSC:} 65P10, 65L05

\end{abstract}

\section{Introduction}

In this paper, we are concerned with the numerical energy-preserving
analysis over a long term of extended Runge--Kutta--Nystr\"{o}m
(ERKN) integrators when applied to the following   highly
oscillatory  Hamiltonian system
\begin{equation}
 \left\{\begin{aligned}
 &\dot{q}=\nabla_p H(q,p),\qquad \ \ q(0)=q^{0},\\
 &\dot{p}=-\nabla_q H(q,p),\qquad  p(0)=p^{0}, \end{aligned}\right.
\label{H-s}%
\end{equation}
with the Hamiltonian
\begin{equation}H(q,p)=\frac{1}{2}\big(\norm{p}^2+\norm{\Omega q}^2\big)+U(q),\label{H}%
\end{equation}
where the vectors $p = (p_1, p_2) \in \mathbb{R}^{d_1}\times
\mathbb{R}^{d_2}$ and $q = (q_1, q_2)\in \mathbb{R}^{d_1}\times
\mathbb{R}^{d_2}$ are partitioned  subject to the partition of the
square matrix
\begin{equation}\Omega= \left(
                          \begin{array}{cc}
                            0_{d_1\times d_1}  & 0 _{d_2\times d_2}  \\
                            0_{d_1\times d_1} & \omega I_{d_2\times d_2} \\
                          \end{array}
                        \right) ,
\label{Omega}%
\end{equation}
and  $\omega$ is a large positive parameter. It is assumed in this
paper that the initial values of \eqref{H-s} satisfy  the condition
\begin{equation}
\frac{1}{2} \norm{p(0)}^2+\frac{1}{2} \norm{\Omega q(0)}^2\leq E, \label{energy condition}%
\end{equation}
where the constant $E$ is independent of $\omega$.  As is known, the
oscillatory energy of the system \eqref{H-s} is
\begin{equation}
I(q,p)=\frac{1}{2}p_2^\intercal p_2 +\frac{1}{2}\omega^2q_2^\intercal q_2. \label{oscillatory energy}%
\end{equation}

These highly oscillatory  Hamiltonian systems  frequently  arise in
a wide variety of applications including applied mathematics,
quantum physics, classical mechanics, molecular biology, chemistry,
astronomy, and  electronics (see, e.g.
\cite{hairer2006,Hochbruck2010,wubook2018,wu2013-book}). Over the
past two decades (and earlier), some novel approaches have been
extensively studied for the highly oscillatory Hamiltonian system
and  we refer the reader to
\cite{franco2006,Grimm2006,Hairer09,Hochbruck1999,Li_Wu(sci2016),Li_Wu(na2016),wang-2016,wu2017-JCAM,wang2017-JCM}
as well as the references contained therein.

The numerical long-time near-conservation of energy for such
equations has already been researched for various numerical
integrators, such as  for symmetric trigonometric integrators in
\cite{Cohen15,Cohen03,Hairer00,hairer2006}, for the
St\"{o}rmer--Verlet method in \cite{Hairer00-1},   for an
implicit-explicit method in \cite{McLachlan14,Stern09},  and for
heterogeneous multiscale methods in \cite{Sanz-Serna09}. The
essential technology  used  in the analysis is  the modulated
Fourier expansion, which was firstly developed  in \cite{Hairer00}
to show long-time almost conservation properties  of numerical
methods for highly oscillatory Hamiltonian systems. A similar
modulated Fourier expansion  was independently presented in
\cite{Guzzo01} on the spectral formulation of the Nekhoroshev
theorem for quasi-integrable Hamiltonian systems. Following those
pioneering  researches, modulated Fourier expansions have been
developed as an important mathematical tool in studying the
long-time   behaviour of numerical methods and differential
equations (see, e.g.
\cite{Cohen06,Sanz-Serna99,Gaucklerthesis,Hairer09,Hairer12-1,Hairer13,iserles08}).
 With  this technology, much work has been done  about
different numerical methods for various systems.  Besides the work
stated above for the highly oscillatory Hamiltonian system
\eqref{H},  the numerical energy conservation for other kinds of
systems has also been researched in many publications, such as for
multi-frequency  Hamiltonian systems in \cite{Cohen05,Cohen06-1},
wave equations in \cite{Cohen08,Cohen08-1,Gauckler16,Hairer08},
Schr\"{o}dinger equations in \cite{Cohen12,Gauckler10-1,Gauckler10},
highly oscillatory Hamiltonian systems without   any non-resonance
condition in \cite{Gauckler13}, and  Hamiltonian systems with a
solution-dependent high frequency in \cite{Hairer16}.

On the other hand,    in order to  effectively solve \eqref{H} in
the sense of structure-preservation,  the authors in \cite{wu2010-1}
formulated a standard form of  extended Runge--Kutta--Nystr\"{o}m
(ERKN) integrators and derived the corresponding order conditions by
the B-series theory associated with   the extended Nystr\"{o}m
trees. In \cite{wu2013-ANM}, the error bounds for explicit ERKN
integrators were researched. Recently, symplecticity conditions of
ERKN integrators were  derived  in \cite{wu2012} and symmetry
conditions were derived in \cite{wu2013-book}. Based on these
conditions, some practical symmetric or/and symplectic ERKN
integrators were constructed and analyzed in \cite{wang2017-Cal}.
The results of numerical experiments appearing in the work mentioned
above have shown that the symmetric or/and symplectic   ERKN
integrators behave very well even in a long-time interval. However,
the theoretical analysis of energy behaviour over  a long term  of
symmetric or  symplectic ERKN integrators has not been considered
and researched yet in the literature, which motives this paper.


The  main  contributions of this work are to show  the long-time
energy behaviour  not only for one-stage explicit symmetric ERKN
integrators but also for symplectic ERKN integrators and to derive
modulated Fourier expansions for non-symmetric methods. Similar
results have been obtained for symmetric trigonometric integrators
in \cite{Cohen15,Cohen05,Hairer00,hairer2006}. However, in this
paper we prove the long-time result  for more diverse methods than
that for those considered previously. In particular, we present the
analysis for both symmetric and symplectic ERKN integrators. We
neither assume symplecticity for symmetric methods, nor assume
symmetry for symplectic methods. It follows from the analysis that
both symmetry and symplecticity can produce a good long-time energy
conservation for ERKN integrators, which means that symmetry and
symplecticity play a similar role in the numerical energy-preserving
behaviour. This is a  new discovery which is of great importance to
geometric integration for highly oscillatory Hamiltonian systems.
Moreover, in contrast to \cite{Cohen15,Cohen05,Hairer00,hairer2006},
in the long-term analysis of symplectic integrators, the formulation
of modulated Fourier expansions does not rely on the symmetry of the
methods. This is of major importance in the context of long-term
analysis of non-symmetric methods. It is also a main conceptual
difference in comparison with
\cite{Cohen15,Cohen05,Hairer00,hairer2006}.

For the analysis of symmetric integrators, we have noted that some
ERKN integrators can be formulated as a Strang splitting method
applied to an averaged equation (see, \cite{Blanes2015}). Very
recently, the authors in \cite{Buchholz2018}  proved second-order
error bounds of trigonometric integrators on the basis of the
interpretation of trigonometric integrators as splitting methods for
averaged
 equations.   Following this way,  the long-term
 analysis of symmetric ERKN integrators will be proved   concisely by exploring the
connection between symmetric ERKN integrators and trigonometric
integrators researched in
\cite{Cohen15,Cohen05,Hairer00,hairer2006}. However, in the analysis
of symplectic ERKN integrators, we do not require the symmetry of
methods and  it is known that  the symmetry  plays an important role
in the construction of modulated Fourier expansions and long-term
analysis given in \cite{Cohen15,Cohen05,Hairer00,hairer2006}.
Therefore, unfortunately, the
  approach used for symmetric ERKN integrators  does not apply to  symplectic ERKN integrators any more.
In order to overcome this difficulty, we  will use the technology of
modulated Fourier expansion developed by Hairer and Lubich  in
\cite{Hairer00} with some novel adaptations for non-symmetric
methods. The modulated Fourier expansion of symplectic ERKN
integrators will be derived
 and two almost-invariants will be
shown. Then  the long-term result can be obtained.

This paper is organized as follows. We first present some
preparatories of ERKN integrators in Section \ref{sec:Formulation}.
The main results as well as  an  illustrative numerical experiment
are given in Section   \ref{sec:examples}. Then in Section \ref{sec:
proof}, the result for symmetric integrators is proved by exploring
the connection between symmetric ERKN integrators and trigonometric
integrators and by using the previous results of symmetric
trigonometric integrators shown in
\cite{Cohen15,Cohen05,Hairer00,hairer2006}. Section \ref{sec:
symplectic proof} gives the proof of long-term result for symplectic
ERKN integrators, where the modulated Fourier expansion  is
constructed for symplectic integrators and two almost-invariants of
the modulated Fourier expansions are studied. The concluding remarks
are made in the last section.

\section{ERKN integrators} \label{sec:Formulation}
The highly oscillatory Hamiltonian system \eqref{H} can be rewritten
as the following  system of  second-order differential equations
\begin{equation}
q^{\prime\prime}(t)+\Omega^2q(t)=g(q(t)),  \qquad
q(0)=q^0,\ \ q'(0)=p^0,\label{prob}%
\end{equation}
where $g$ is the negative gradient of  a  real-valued function
$U(q)$.  ERKN integrators were  first formulated for integrating
  (\ref{prob}) in \cite{wu2010-1} and here we  summarize  the
scheme of one-stage  ERKN  integrators  as follows.

\begin{mydef}
\label{erkn}  (See \cite{wu2010-1}) A one-stage  explicit ERKN
integrator for solving (\ref{prob}) is defined by%
 \begin{equation}
\begin{array}
[c]{ll}%
Q^{n+c_{1}} &
=\phi_{0}(c_{1}^{2}V)q^{n}+hc_{1}\phi_{1}(c_{1}^{2}V)p^{n},\\
q^{n+1} & =\phi_{0}(V)q^{n}+h\phi_{1}(V)p^{n}+h^{2}
 \bar{b}_{1}(V)g(Q^{n+c_{1}}),\\
p^{n+1} & =-h\Omega^2\phi_{1}(V)q^{n}+\phi_{0}(V)p^{n}+h\textstyle
b_{1}(V)g(Q^{n+c_{1}}),
\end{array}
  \label{erknSchems}%
\end{equation}
where   $h$ is  a  stepsize, $c_1\in[0,1]$ is a real constant,
$b_{1}(V)$ and $\bar{b}_{1}(V)$  are matrix-valued and bounded
functions of $V\equiv h^{2}\Omega^2$, and
\begin{equation}\label{Phi-j}
\phi_{j}(V):=\sum\limits_{k=0}^{\infty}\dfrac{(-1)^{k}V^{k}}{(2k+j)!},\qquad j=0,1.%
\end{equation}
\end{mydef}

From \eqref{Phi-j}, it is clear that
\begin{equation*}
\phi_{0}(V)=\cos(\sqrt{V})=\cos(h\Omega),\qquad  \phi_{1}(V)=
\sin(\sqrt{V}) (\sqrt{V})^{-1} = \sinc(h\Omega),
\end{equation*}
where $\sinc(h\Omega) =\sin(h\Omega)(h\Omega)^{-1}.$ Thence  the
scheme of one-stage explicit   ERKN integrators for \eqref{prob} can
be reformulated as follows.
\begin{mydef}
\label{numerical method}  The one-stage explicit ERKN integrator
 for integrating  \eqref{prob}   is given by
\begin{equation}
\begin{aligned} Q^{n+c_{1}} &
=\cos(c_1h\Omega)q^{n}+c_{1}h\sinc(c_1h\Omega)p^{n},\\
q^{n+1} & =\cos(h\Omega)q^{n}+h\sinc(h\Omega)p^{n}+h^{2}
 \bar{b}_{1}(h\Omega)g(Q^{n+c_{1}}),\\
p^{n+1} &
=-h\Omega^2\sinc(h\Omega)q^{n}+\cos(h\Omega)p^{n}+h\textstyle
b_{1}(h\Omega)g(Q^{n+c_{1}}),
\end{aligned}
\label{methods}%
\end{equation}
where the functions $\bar{b}_1(h\Omega)$ and $b_1(h\Omega)$ are
real-valued and bounded functions of $h\Omega$.


\end{mydef}

 As shown  in \cite{wu2012,wu2013-book}, we obtain the following
conditions for the integrator \eqref{methods}  to be symmetric and
symplectic.

\begin{mytheo}\label{symmetric thm}
The ERKN integrator \eqref{methods} is symmetric if and only if
\begin{equation}\begin{aligned}\label{sym cond}&c_1=1/2,\quad   (I+\cos(h\Omega))\bar{b}_1(h\Omega)=\mathrm{sinc}(h\Omega)b_1(h\Omega) .
\end{aligned}\end{equation}
\end{mytheo}
\textbf{Proof.} It follows from the symmetry conditions given in
\cite{wu2013-book} that this integrator is symmetric if and only if
\begin{equation}\begin{aligned}\label{1s2con1}&c_1=1/2,\
\bar{b}_{1}(V)=\phi_{1}(V)b_{1}(V)-\phi_{0}(V)\bar{b}_{1}(V),\\
&\phi_{0}(c_{1}^{2}V)\bar{b}_{1}(V)=c_{1}\phi_{1}(c_{1}^{2}V)b_{1}(V).
\end{aligned}\end{equation}
By solving the second equation in \eqref{1s2con1}, we obtain  the
second result of \eqref{sym cond}. It  can also be  verified that
under the condition  \eqref{sym cond}, the third equation of
\eqref{1s2con1} is true.

\begin{mytheo}\label{symplectic thm}
 For any  real number $d_1$, if the coefficients are
determined by
\begin{equation}\begin{aligned}\label{symple cond}&
b_1(h\Omega)=d_1\cos((1-c_1)h\Omega),\
 \ \ \bar{b}_1(h\Omega)=d_1(1-c_1)\mathrm{sinc}((1-c_1)h\Omega),
\end{aligned}\end{equation}
then the ERKN integrator \eqref{methods} is symplectic.
\end{mytheo}
\textbf{Proof.} According to the symplectic conditions  given in
\cite{wu2012}, we know  that this method is symplectic if  the
following equations
\begin{equation*}
\begin{array}
[c]{ll}%
\phi_{0}(V)b_{1}(V)+V\phi_{1}(V)\bar{b}_{1}(V)=d_{1}\phi_{0}(c_{1}^{2}V),\ \ \ d_1\in \mathbb{R},\\
\phi_{1}(V)b_{1}(V)-\phi_{0}(V)\bar{b}_{1}(V)=c_{1}d_{1}\phi_{1}(c_{1}^{2}V),
\end{array}
\end{equation*}
are satisfied.  The result is  directly  obtained by solving these
two equations.

\begin{myrem}
These two theorems confirm the fact that an  ERKN integrator can be
symmetric and symplectic, or symmetric but not symplectic, or
symplectic but not symmetric.
\end{myrem}

As some examples of ERKN integrators with certain structure
characteristics, we present six practical one-stage explicit
integrators and their coefficients are
 listed  in Table \ref{praERKN}. By Theorems \ref{symmetric thm} and
\ref{symplectic thm}, it can be verified that ERKN1 is neither
symmetric nor symplectic,   ERKN2 is symmetric and symplectic,
ERKN3-4    are symmetric but not symplectic, and ERKN5-6    are
symplectic but not symmetric.
\renewcommand\arraystretch{1.8}
\begin{table}[!htb]$$
\begin{array}{|c|c|c|c|c|c|c|c|}
\hline
\text{Methods} &c_1  &\bar{b}_1(h\omega)   &b_1(h\omega)  & \text{Symmetric}  &  \text{Symplectic}  \\
\hline
\text{ERKN1} & \frac{1}{2} & \frac{1}{2} \textmd{sinc}^2(\frac{h\omega}{2}) & \cos(\frac{h\omega}{2})   & \text{Non}  &\text{Non}   \cr
 \text{ERKN1} & \frac{1}{2} &  \frac{1}{2} \textmd{sinc}(\frac{h\omega}{2})  & \cos(\frac{h\omega}{2})    & \text{Symmetric} & \text{Symplectic} \cr
 \text{ERKN3}   & \frac{1}{2} &\frac{1}{2}\textmd{sinc}(h\omega)
\cos(\frac{h\omega}{2}) & \cos^3(\frac{h\omega}{2})
& \text{Symmetric} & \text{Non} \cr
\text{ERKN4} & \frac{1}{2} &  \frac{1}{2} \textmd{sinc}^2(\frac{h\omega}{2})  &\textmd{sinc}(\frac{h\omega}{2}) \cos(\frac{h\omega}{2})   & \text{Symmetric} & \text{Non} \cr
\text{ERKN5} & \frac{2}{5} &  \frac{3}{5} \textmd{sinc}(\frac{3h\omega}{5})   & \cos(\frac{3h\omega}{5})    &\text{Non} & \text{Symplectic} \cr
\text{ERKN6} & \frac{1}{5} &  \frac{4}{5} \textmd{sinc}(\frac{4h\omega}{5})   &\cos(\frac{4h\omega}{5})   & \text{Non}& \text{Symplectic} \cr
 \hline
\end{array}
$$
\caption{Six one-stage explicit ERKN integrators.} \label{praERKN}
\end{table}

\section{Main results and numerical examples}\label{sec:examples}
Before presenting the main results of this paper, we make the
following assumptions.
\begin{assum}\label{ass}
\begin{itemize}
\item  Assume that the initial values satisfy \eqref{energy
condition}.

\item The numerical solution is assumed to stay in a compact set.

\item A lower bound on the  stepsize  is posed as:
\begin{equation}
h\omega \geq c_0 > 0. \label{h condition}%
\end{equation}
\item The numerical non-resonance condition is assumed to be
 held
\begin{equation}
|\sin(\frac{1}{2}kh\omega)| \geq c \sqrt{h}\ \ \mathrm{for} \ \ k=1,2,\ldots,N\ \   \mathrm{with} \ \ N\geq2.\label{numerical non-resonance cond}%
\end{equation}
For a given $h$ and $\omega$, this condition imposes a restriction
on $N$. In the following, $N$ is a fixed integer such that
\eqref{numerical non-resonance cond} holds.

\item For the coefficients of the ERKN integrators, it is  assumed that the function
\begin{equation}
 \sigma (\nu)=\frac{\mathrm{sinc}(\nu)\cos(\frac{1}{2}\nu)}{2\bar{b}_1(\nu)}+\nu^2\mathrm{sinc}(
\nu)\frac{\frac{1}{2} \mathrm{sinc}(\frac{1}{2}\nu)
}{2b_1(\nu)}\label{sigma}%
\end{equation}
is bounded from below and above:
\begin{equation}0<c_1\leq \sigma (\nu) \leq C_1\quad \mathrm{for} \quad \nu=0, h\omega,\label{sigma bound}%
\end{equation}
or the same estimate holds for $-\sigma$ instead of $\sigma$.

\end{itemize}
\end{assum}

It is noted that the first four assumptions are considered by many
publications in the  energy analysis of symmetric trigonometric
integrators for the Hamiltonian  system \eqref{H-s} (see, e.g.
\cite{Cohen06,Cohen15,Hairer00,hairer2006}). The last assumption is
obtained in the remainder analysis of this paper and it is similar
to Assumption B proposed in \cite{Cohen15}.

 With regard to the long-time total and oscillatory
energy conservation along symmetric or symplectic ERKN integrators,
we have the following two main results of this paper, which will be
proved in detail in the next two sections, respectively.
\begin{mytheo} \label{thm: symetric Long-time
thm}  Under the conditions given in Assumption \ref{ass} and   the
symmetry condition \eqref{sym cond}, for one-stage explicit
symmetric ERKN integrators,  it holds that
\begin{equation*}
\begin{aligned}
H(q^n,p^n)&=H(q^0,p^0)+\mathcal{O}(h),\\
I(q^n,p^n)&=I(q^0,p^0)+\mathcal{O}(h)\\
\end{aligned}
\end{equation*}
for $0\leq nh\leq h^{-N+1}.$ The constants symbolized by
$\mathcal{O}$ depend on $N, T$ and the constants in the assumptions,
but are independent of $n,\ h,\ \omega$.
\end{mytheo}

\begin{mytheo} \label{thm: symplec Long-time
thm}  Under the conditions of Assumption \ref{ass} without the
requirement \eqref{sigma bound} and   the symplecticity condition
\eqref{symple cond} with $d_1=1$, for one-stage explicit symplectic
ERKN integrators,   we have
\begin{equation*}
\begin{aligned}
H(q^n,p^n)&=H(q^0,p^0)+\mathcal{O}(h),\\
I(q^n,p^n)&=I(q^0,p^0)+\mathcal{O}(h),
\end{aligned}
\end{equation*}
where $0\leq nh\leq h^{-N+1}.$ The constants symbolized by
$\mathcal{O}$ depend on $N, T$ and the constants in the assumptions,
but are independent of $n,\ h,\ \omega$.
\end{mytheo}

As an  important nonlinear model,  we consider the
Fermi--Pasta--Ulam problem.  This model describes classical and
quantum systems of interacting particles in the physics of nonlinear
phenomena. Denote by $x_{i}$   a scaled displacement of the $i$th
stiff spring and by $x_{m+i}$  a scaled expansion or compression of
the $i$th stiff spring. Their corresponding velocities are expressed
in $y_{i}$ and  $y_{m+i}$, respectively. Then the problem can be
formulated by a Hamiltonian system with the Hamiltonian
 \begin{equation*}
\begin{array}
[c]{ll}%
H(y,x)
=&\dfrac{1}{2}\textstyle\sum\limits_{i=1}^{2m}y_{i}^{2}+\dfrac
{\omega^{2}}{2}\textstyle\sum\limits_{i=1}^{m}x_{m+i}^{2}+\dfrac{1}{4}%
[(x_{1}-x_{m+1})^{4}\\
& +\textstyle\sum\limits_{i=1}^{m-1}(x_{i+1}-x_{m+i-1}-x_{i}-x_{m+i}%
)^{4}+(x_{m}+x_{2m})^{4}].
\end{array}
\end{equation*}

\begin{figure}[ptb]
\centering\tabcolsep=2mm
\begin{tabular}
[l]{lll}%
\includegraphics[width=12cm,height=3cm]{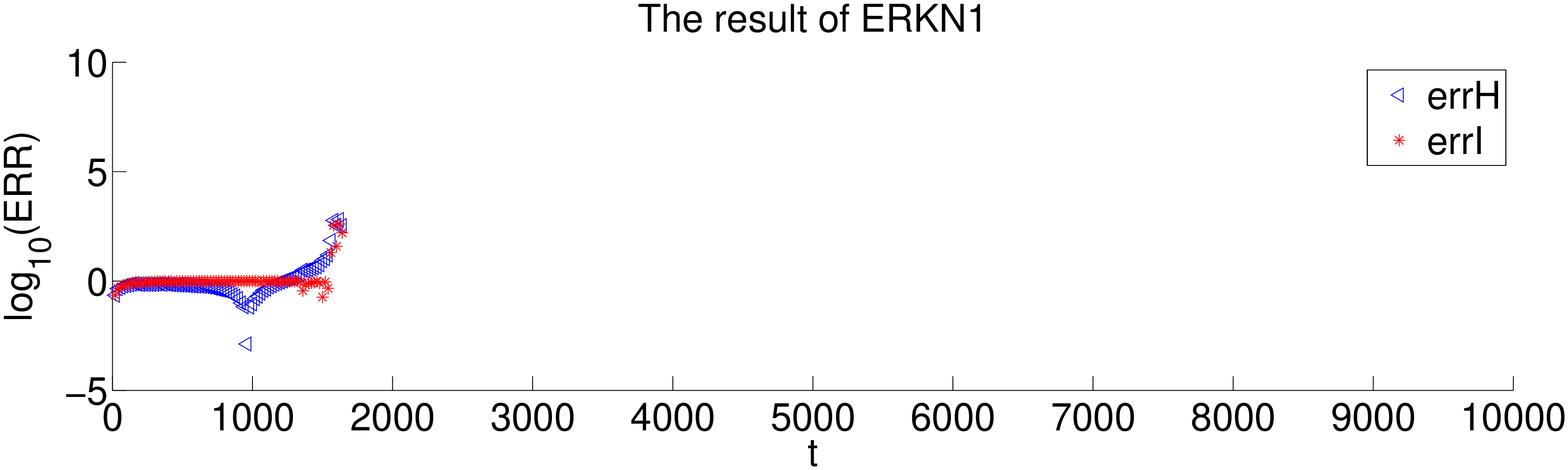}\\
\includegraphics[width=12cm,height=3cm]{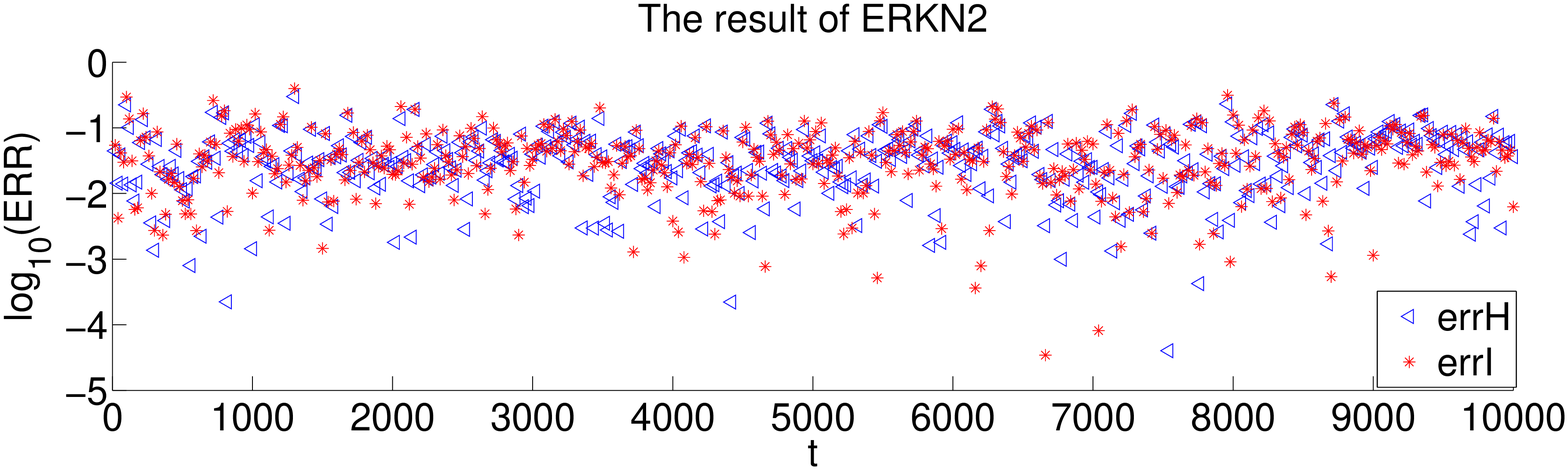}\\
\includegraphics[width=12cm,height=3cm]{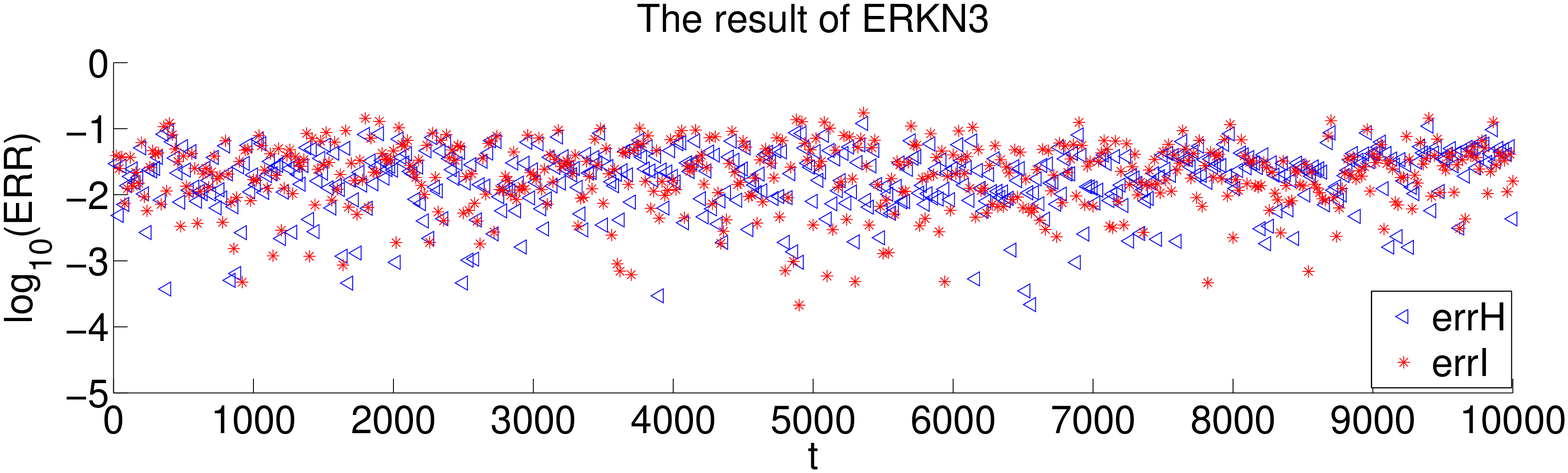}\\
\includegraphics[width=12cm,height=3cm]{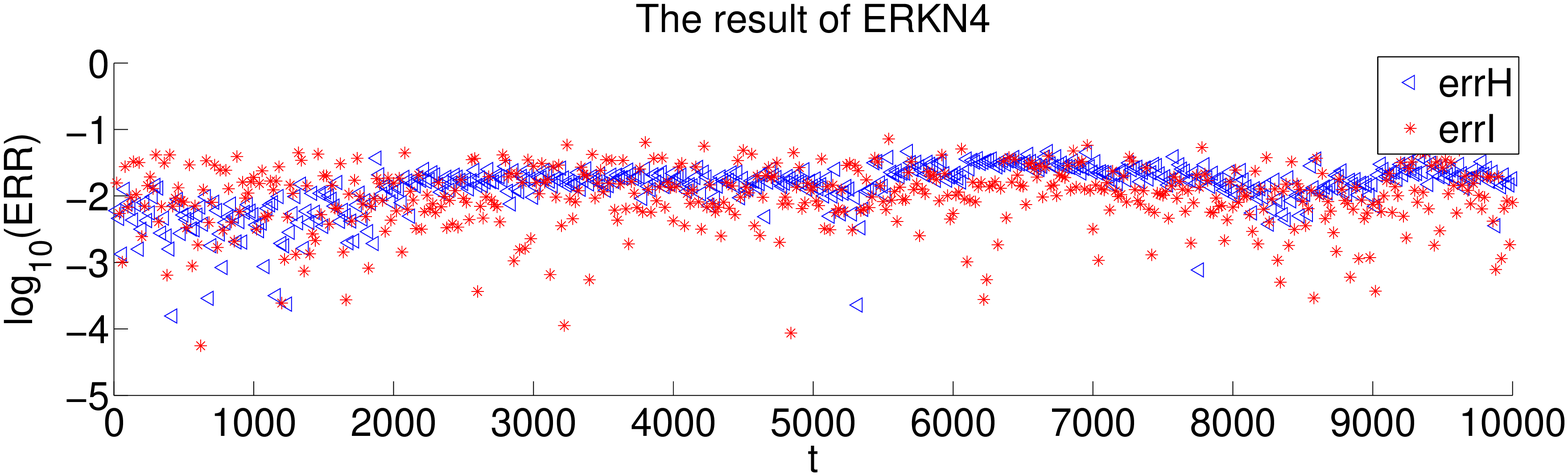}\\
\includegraphics[width=12cm,height=3cm]{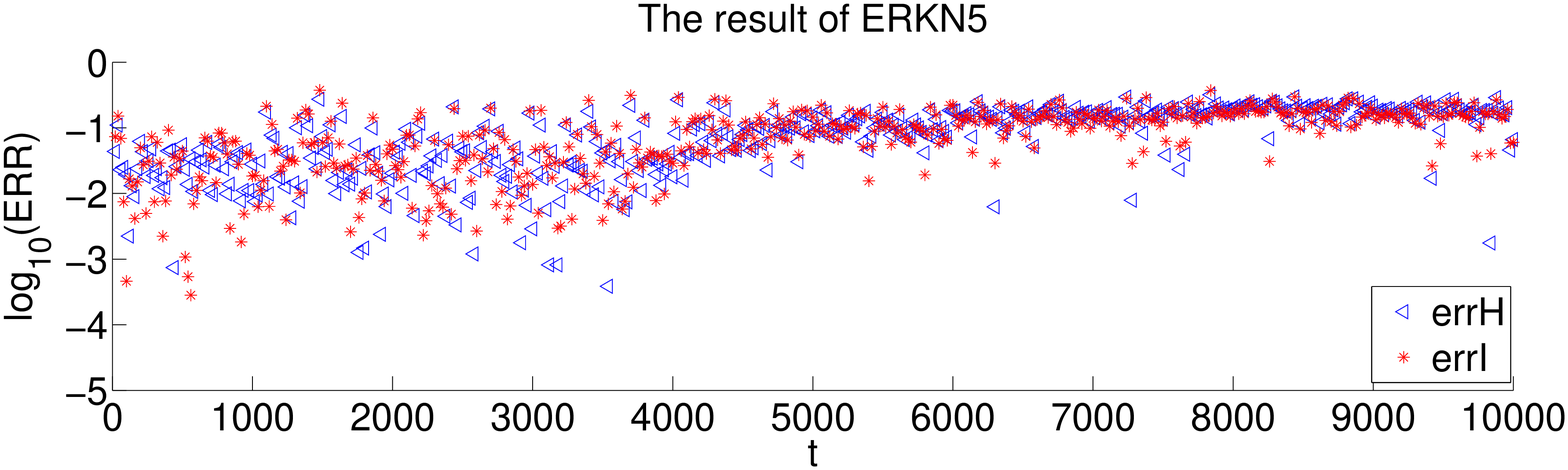}\\
\includegraphics[width=12cm,height=3cm]{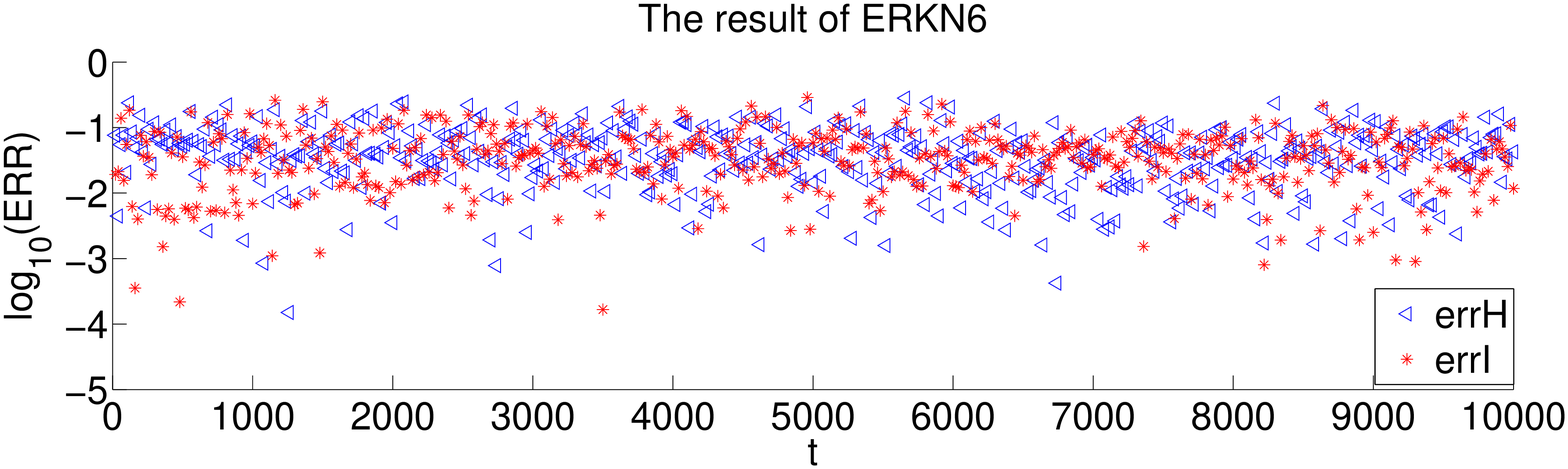}
\end{tabular}
\caption{The logarithm of the  errors ($ERR$)   of  $H$ and   $I$ against $t$ with   $h=0.1$ and $\omega=50$.}%
\label{fig0}%
\end{figure}

\begin{figure}[ptb]
\centering\tabcolsep=2mm
\begin{tabular}
[c]{ccc}%
\includegraphics[width=12cm,height=3cm]{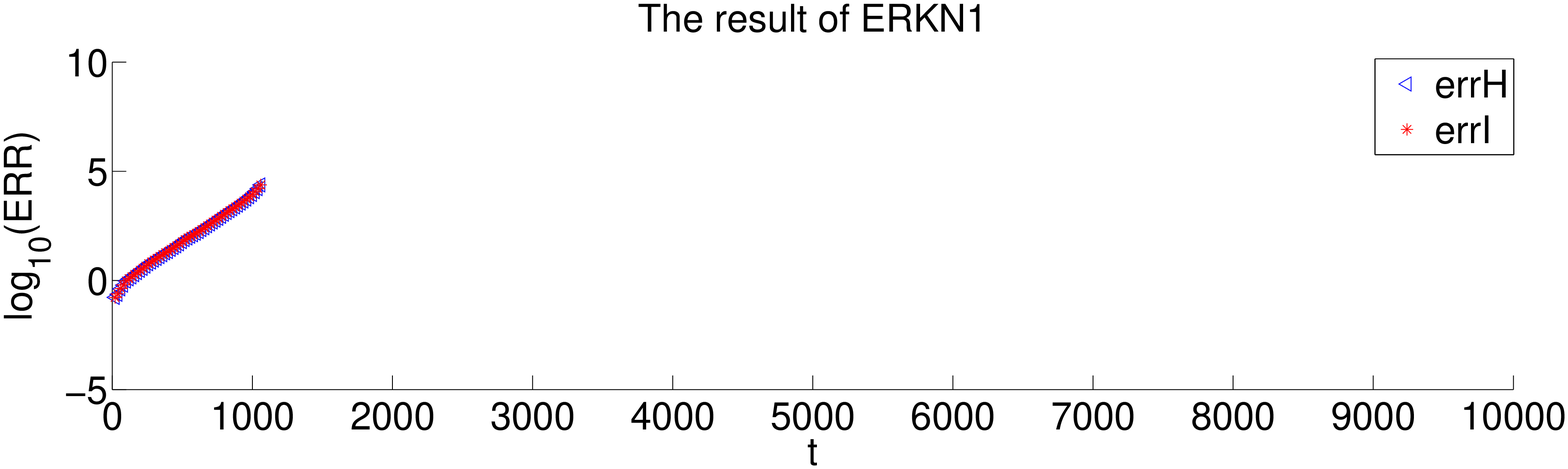}\\
\includegraphics[width=12cm,height=3cm]{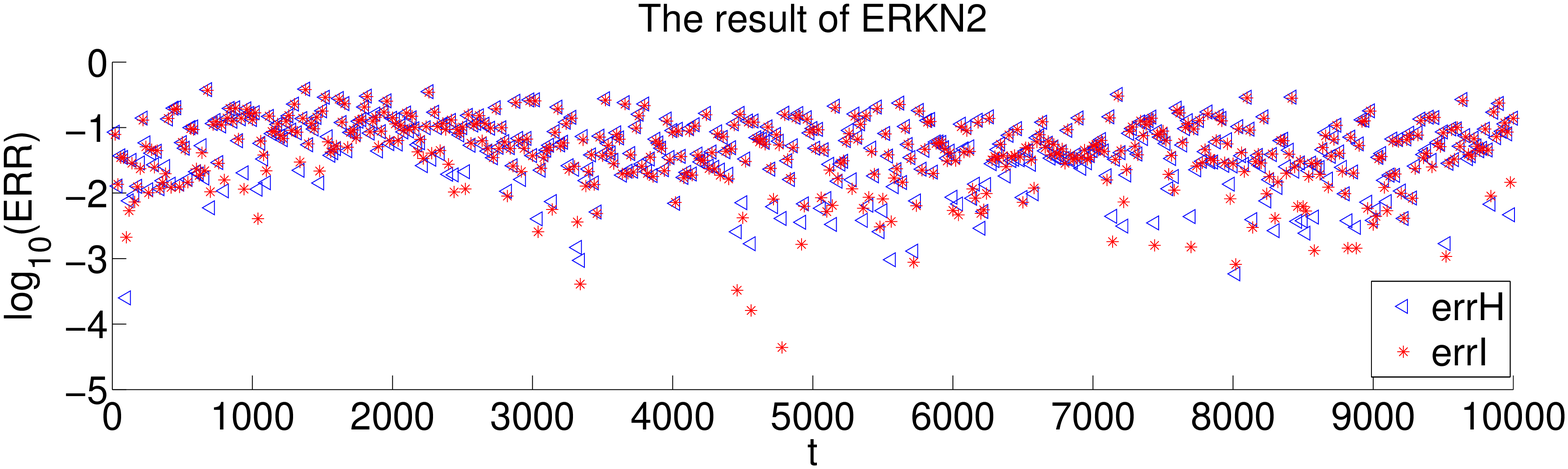}\\
\includegraphics[width=12cm,height=3cm]{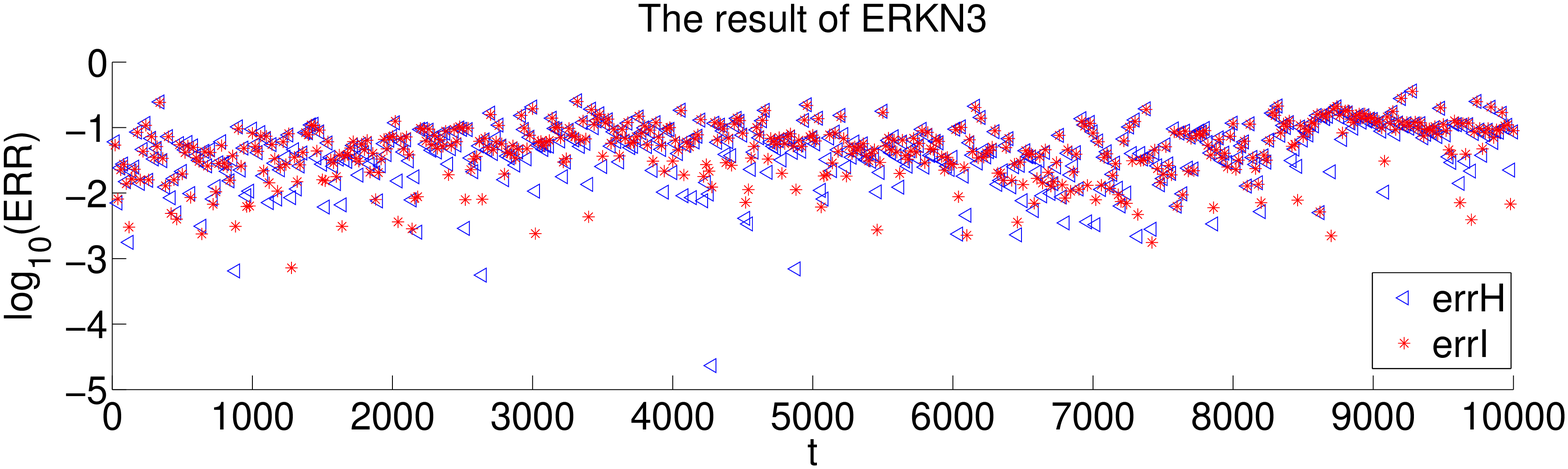}\\
\includegraphics[width=12cm,height=3cm]{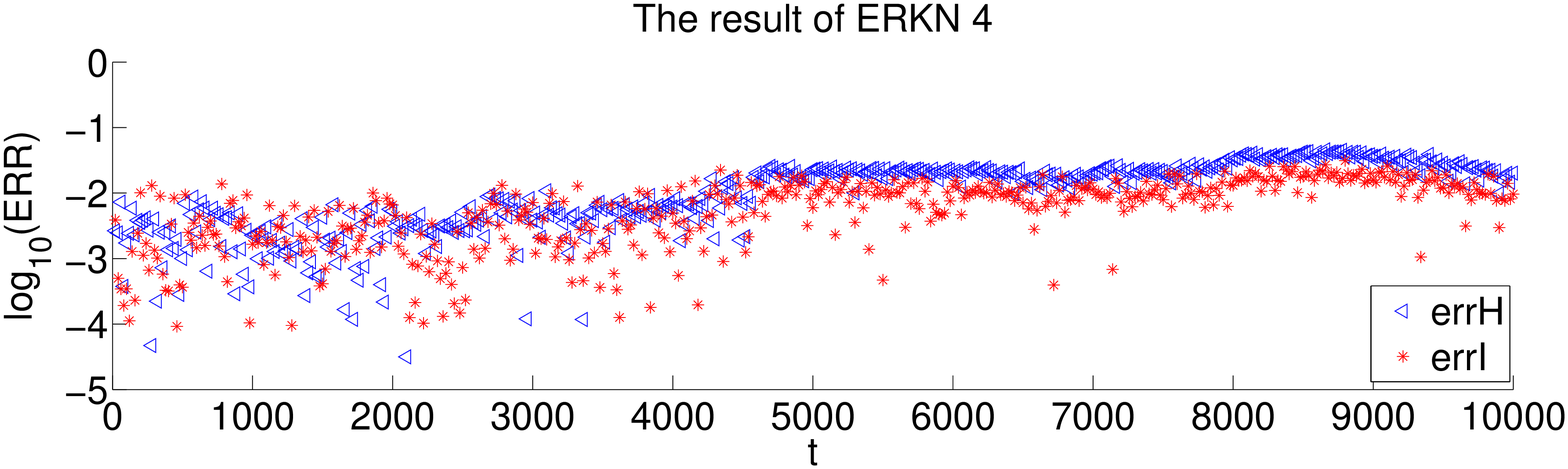}\\
\includegraphics[width=12cm,height=3cm]{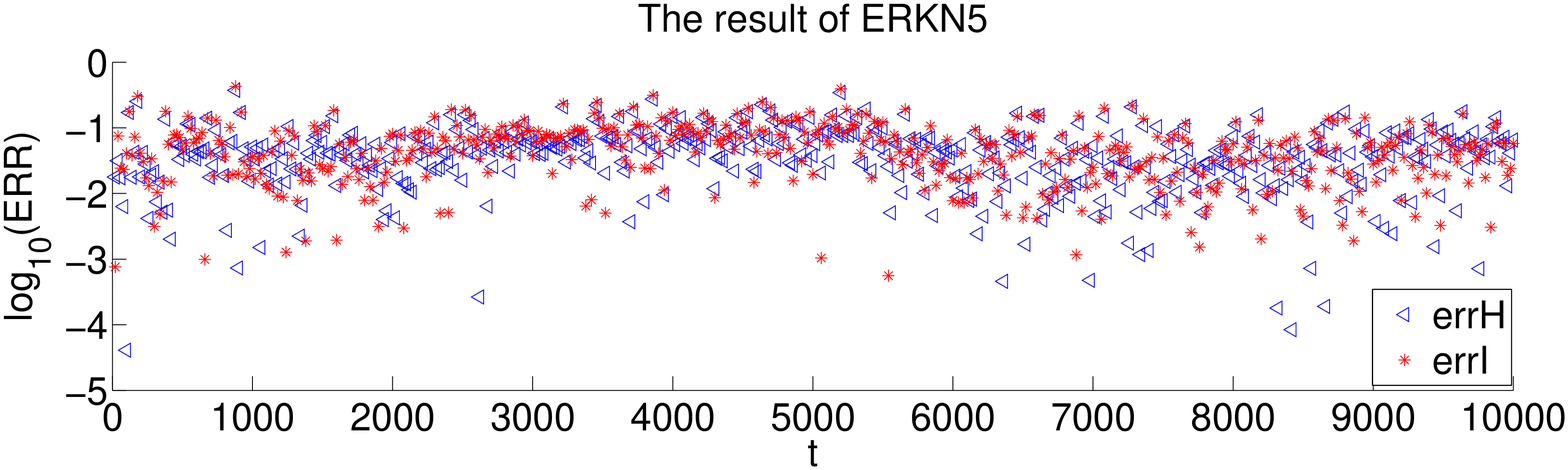}\\
\includegraphics[width=12cm,height=3cm]{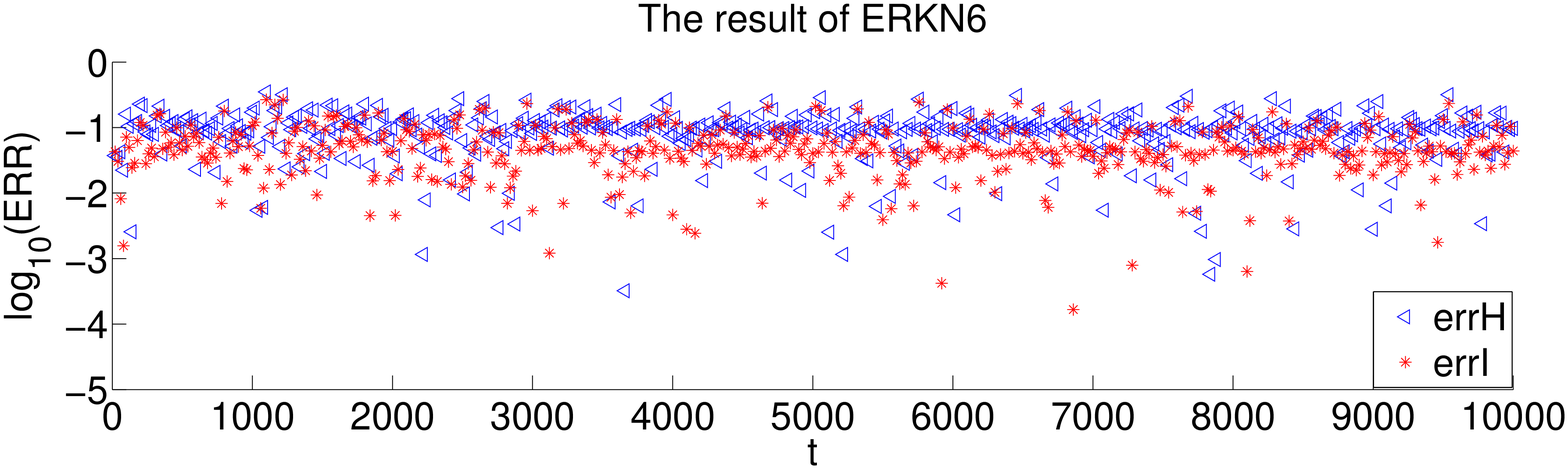}
\end{tabular}
\caption{The logarithm of the  errors ($ERR$)   of  $H$ and   $I$ against $t$ with   $h=0.1$ and $\omega=200$.}%
\label{fig1}%
\end{figure}

\begin{figure}[ptb]
\centering\tabcolsep=2mm
\begin{tabular}
[c]{ccc}%
\includegraphics[width=12cm,height=3cm]{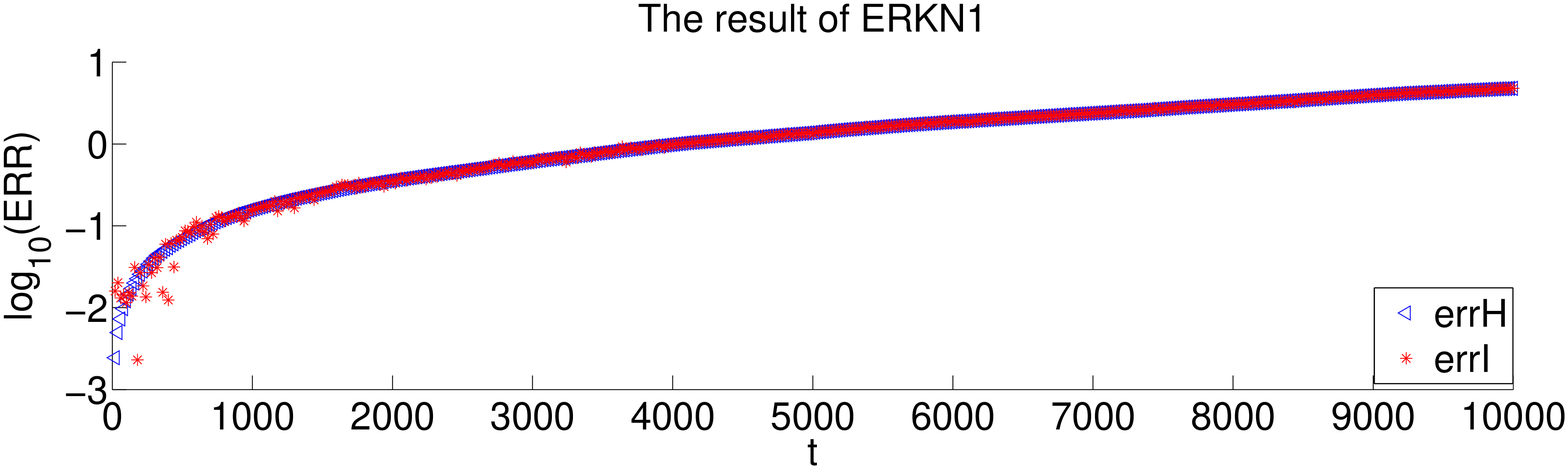}\\
\includegraphics[width=12cm,height=3cm]{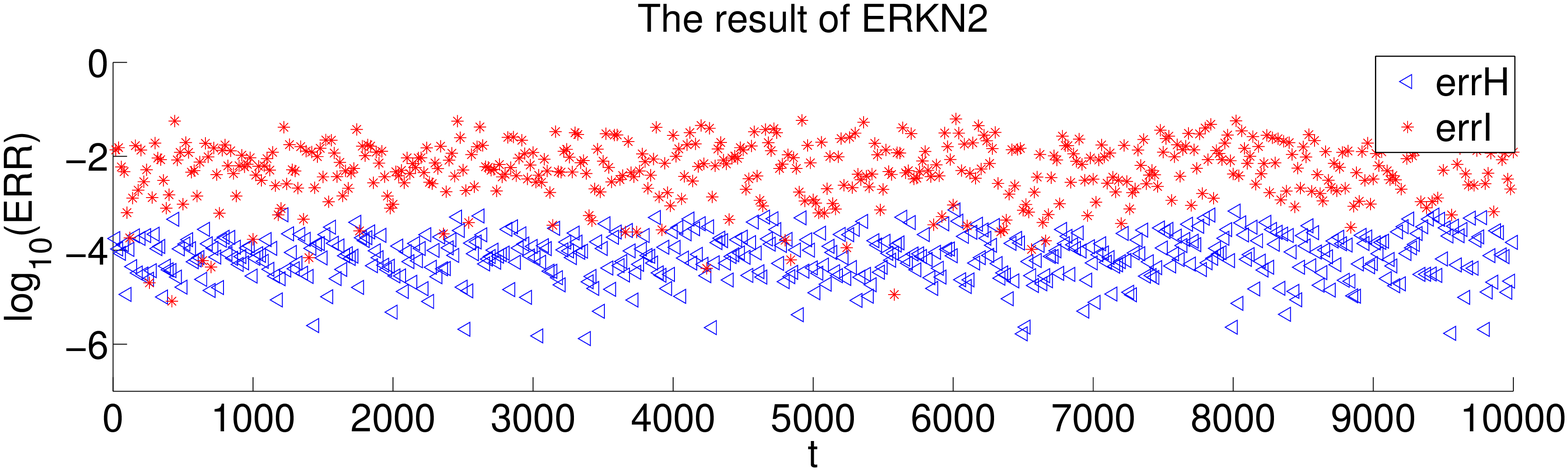}\\
\includegraphics[width=12cm,height=3cm]{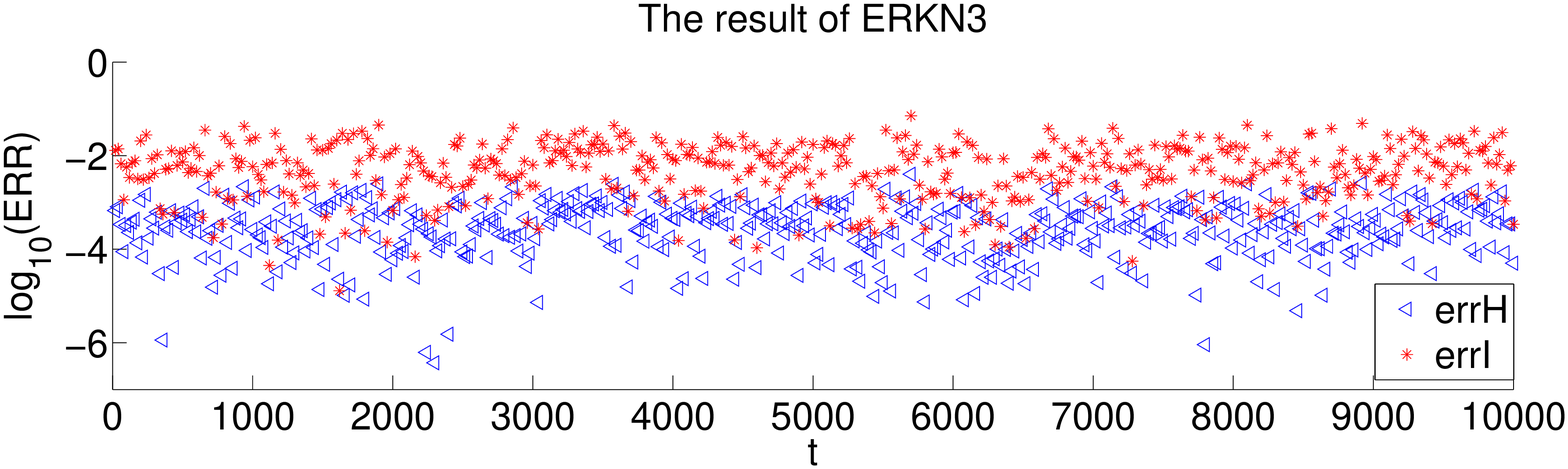}\\
\includegraphics[width=12cm,height=3cm]{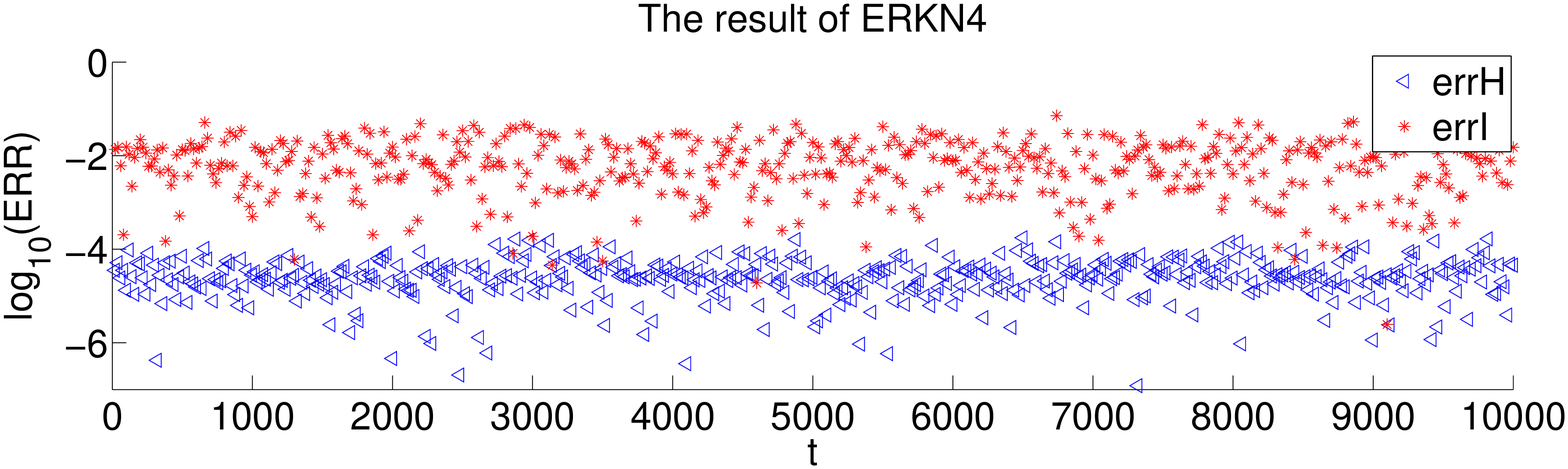}\\
\includegraphics[width=12cm,height=3cm]{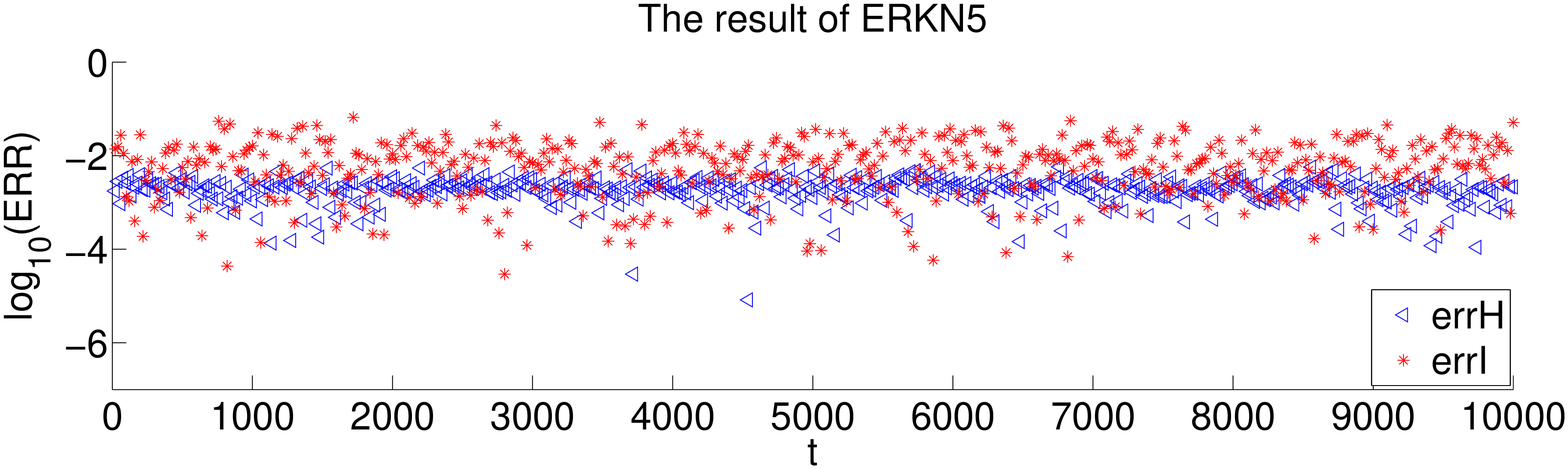}\\
\includegraphics[width=12cm,height=3cm]{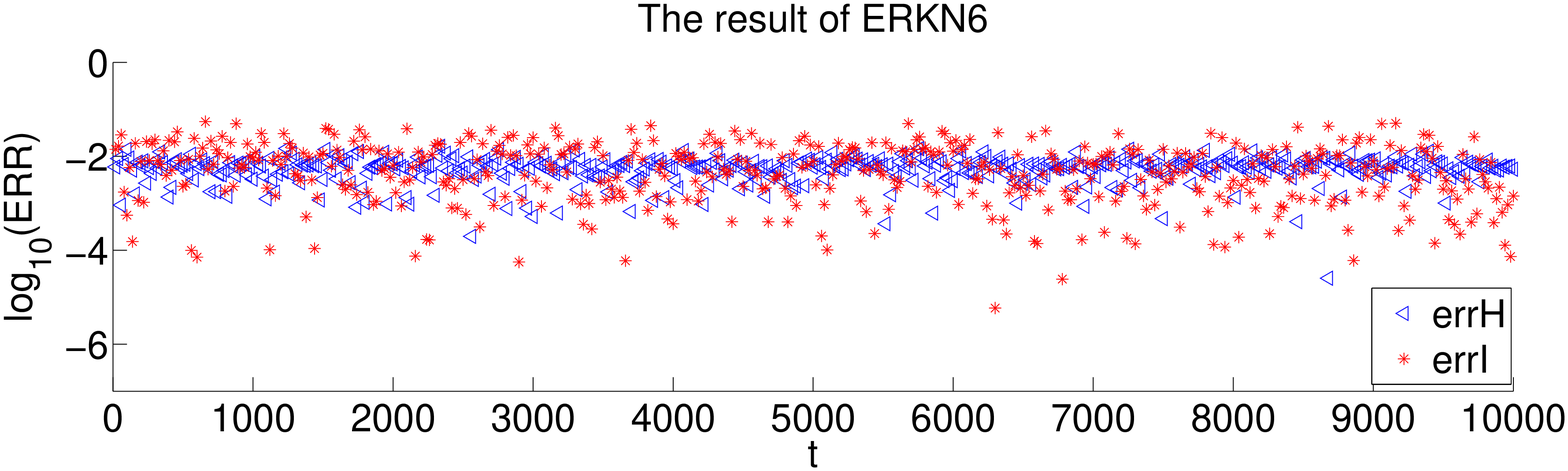}\\
\end{tabular}
\caption{The logarithm of the  errors ($ERR$)   of  $H$ and   $I$ against $t$ with   $h=0.01$ and $\omega=50$.}%
\label{fig1-1}%
\end{figure}

\begin{figure}[ptb]
\centering\tabcolsep=2mm
\begin{tabular}
[c]{ccc}%
\includegraphics[width=12cm,height=3cm]{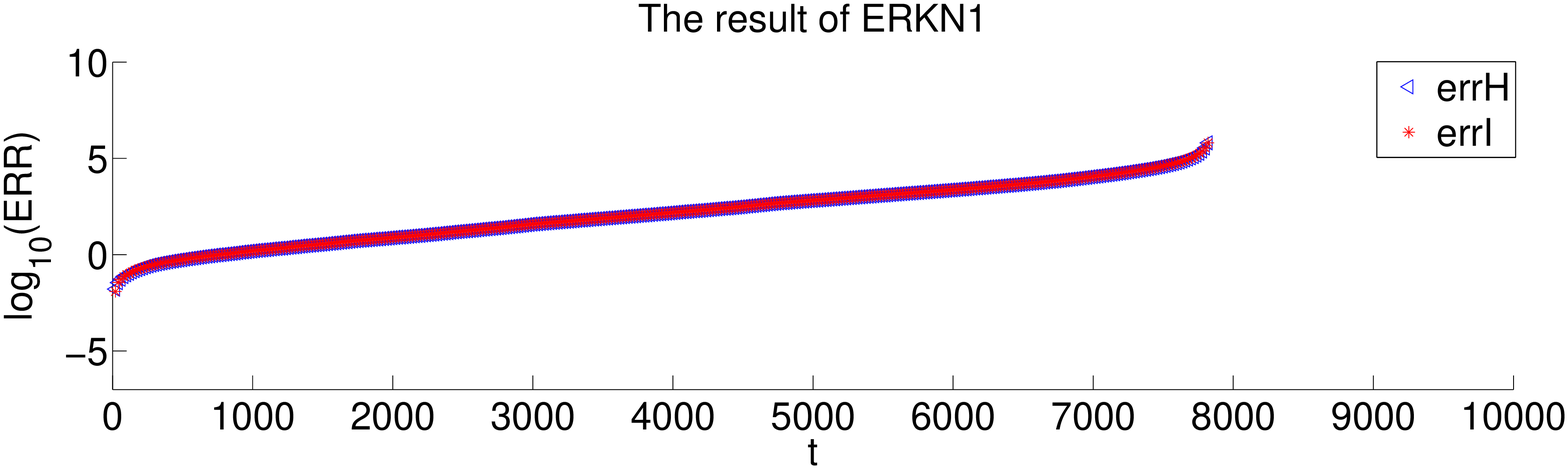}\\
\includegraphics[width=12cm,height=3cm]{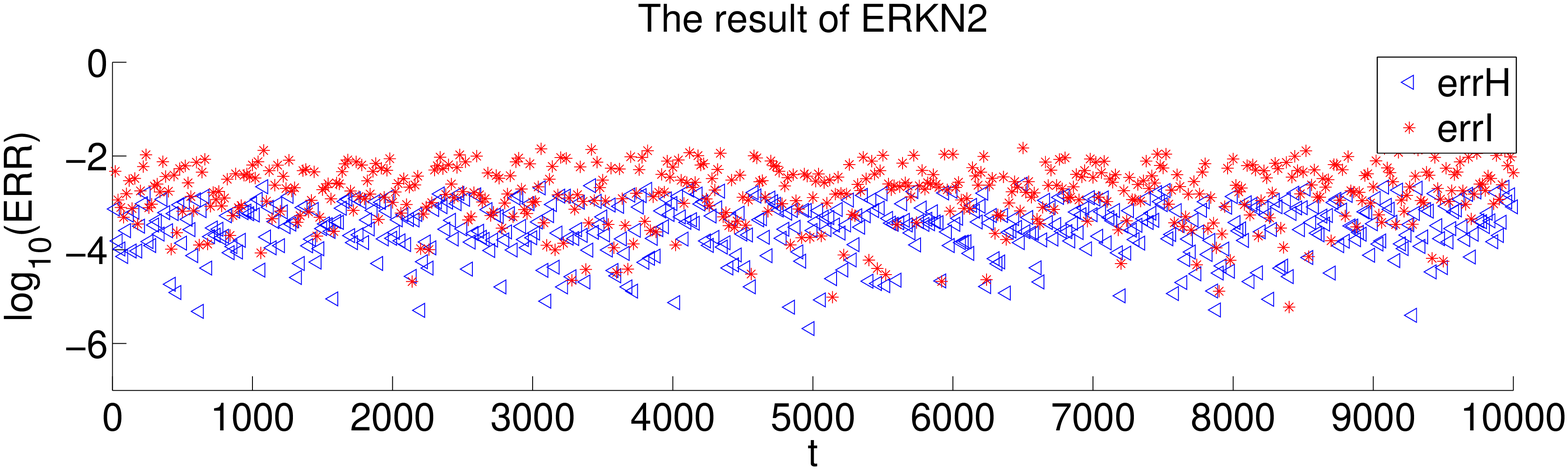}\\
\includegraphics[width=12cm,height=3cm]{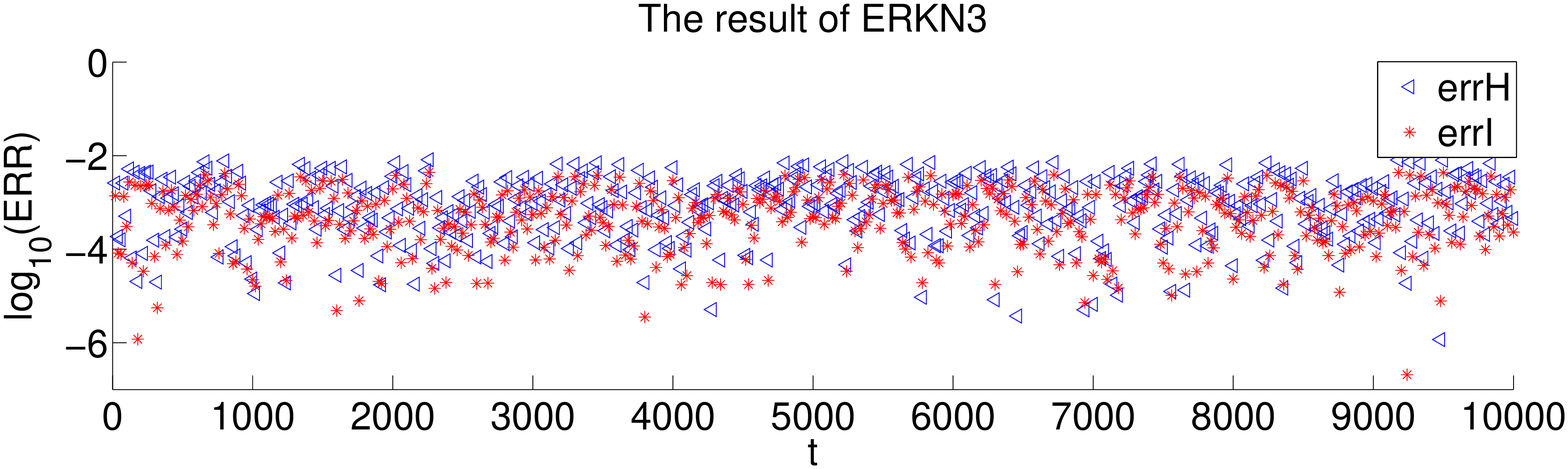}\\
\includegraphics[width=12cm,height=3cm]{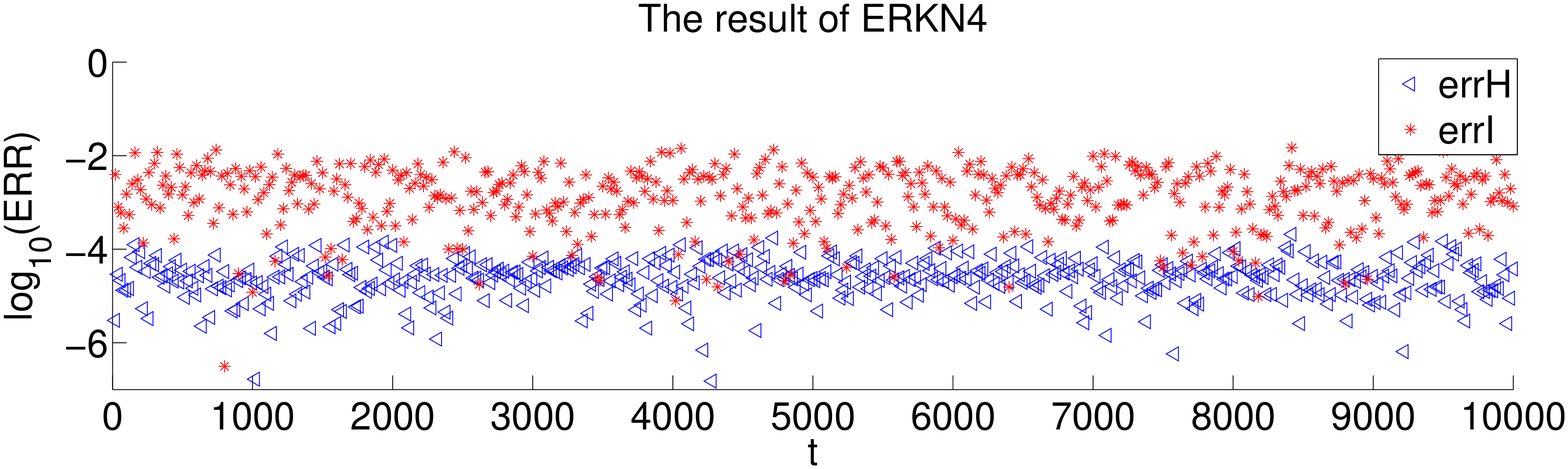}\\
\includegraphics[width=12cm,height=3cm]{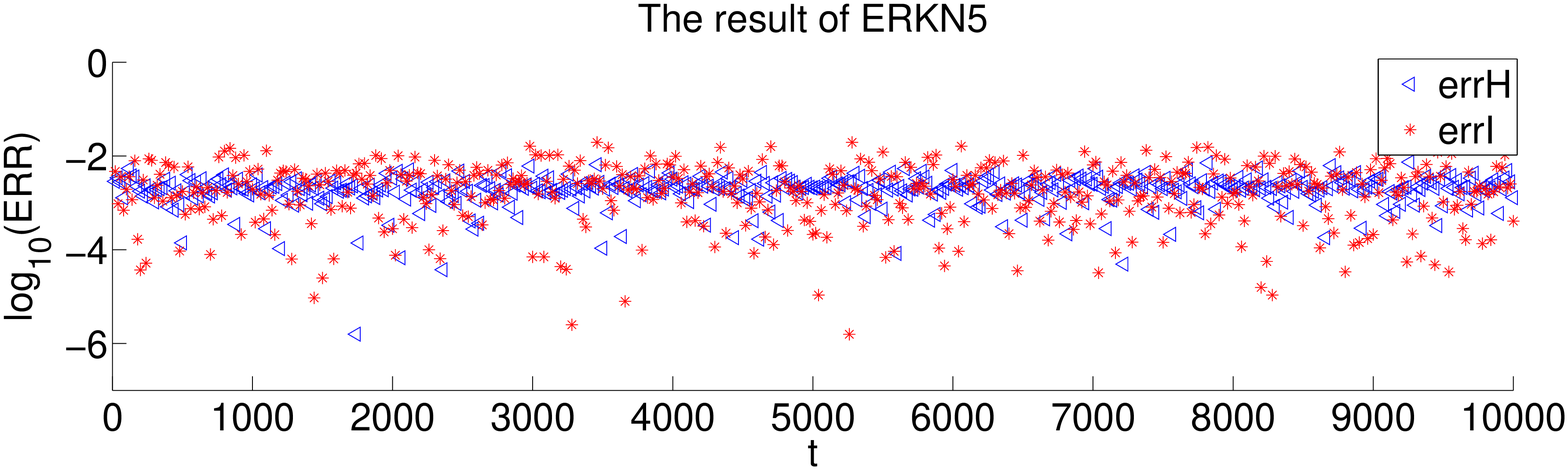}\\
\includegraphics[width=12cm,height=3cm]{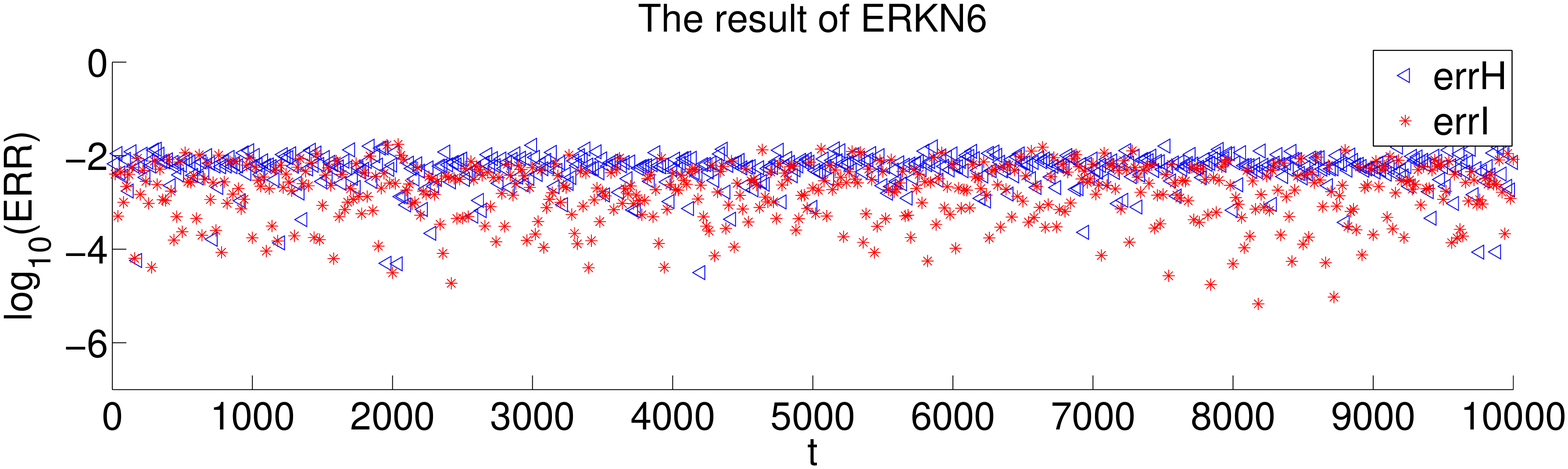}\\
\end{tabular}
\caption{The logarithm of the  errors ($ERR$)   of  $H$ and   $I$ against $t$ with   $h=0.01$ and $\omega=200$.}%
\label{fig2}%
\end{figure}

Following \cite{hairer2006}, we choose $m=3$ and
\begin{equation*}\label{initial date} \ x_{1}(0)=1,\ y_{1}(0)=1,\
x_{4}(0)=\dfrac{1}{\omega},\ y_{4}(0)=1
 \end{equation*}
with zero for the remaining initial values.  First, the system is
integrated in the interval $[0,10000]$ with   $h=0.1$ and
$\omega=50$. The errors of the total energy $H$ and oscillatory
energy $I$ against $t$ for different ERKN integrators are shown in
Fig. \ref{fig0}. Then we increase the size of $\omega$ to $200$ and
see Fig. \ref{fig1} for the results. It can be observed  that the
performance of   ERKN1  is  affected  by $\omega$ and the
corresponding errors become large  as  $\omega$ increases. However,
the other methods are nearly not  affected  by $\omega$. Finally, we
use a smaller stepsize $h=0.01$ and the results are presented in
Figs. \ref{fig1-1}-\ref{fig2}. Under this situation, all the methods
behave better than the former.

It follows from all  the results that the integrators ERKN2-6
approximately conserve $H$ and $I$ quite well over a long term  for
different stepsizes and $\omega$. ERKN1 does not approximately
conserve $H$ and $I$. In the light of  the results, it seems that
the symmetry or symplecticity is   essential for the conservation of
$H$ and $I$. All these numerical  behaviours can be explained   by
the theoretical conclusions. According to Theorems \ref{thm:
symetric Long-time thm} and \ref{thm: symplec Long-time thm},   it
can be concluded that  if ERKN integrators are symmetric or
symplectic, they have a near conservation of $H$ and $I$ over a long
term. This is the theoretical reason for the behaviour shown in this
numerical experiment.


\section{Proof for symmetric integrators} \label{sec: proof}
In this section, we prove the result of long-time total and
oscillatory energy conservation for symmetric ERKN integrators.
\subsection{One important connection}
We denote the numerical flow of a one-stage explicit symmetric ERKN
integrator   \eqref{methods}   with the symmetry condition
\eqref{sym cond} by $\Phi_h$, i.e., $$(q^{n+1},p^{n+1}) = \Phi_h
(q^{n},p^{n}).$$ We now consider a Strang splitting applied to an
averaged equation. More precisely, let $\Phi_{h,\textmd{L}}$ be the
time $h$  flow of the linear equation $\frac{d}{dt}(q,p)=(p,
\Omega^2q)$ and $\Phi_{h,\textmd{NL}}$ be the time $h$  flow of the
nonlinear equation $\frac{d}{dt}(q,p)=(0,\Upsilon(h \Omega)g(q))$.
We then consider the following  Strang splitting method
 \begin{equation*}
\begin{array}[c]{ll}
&\textmd{1. } (q_{+}^{n},p_{+}^{n}) = \Phi_{h/2,\textmd{L}}
(q^{n},p^{n}): \\
&\left(
             \begin{array}{c}
              q_{+}^{n} \\
              p_{+}^{n} \\
             \end{array}
           \right)=\left(
                                                           \begin{array}{cc}
                                                            \cos(\frac{h\Omega}{2})& \Omega^{-1}\sin(\frac{h\Omega}{2}) \\
                                                             -\Omega \sin(\frac{h\Omega}{2}) & \cos(\frac{h\Omega}{2}) \\
                                                           \end{array}
                                                         \right)\left(
             \begin{array}{c}
               q^{n} \\
               p^{n} \\
             \end{array}
           \right),\\
&\textmd{2. }(q_{-}^{n},p_{-}^{n}) = \Phi_{h,\textmd{NL}}
(q_{+}^{n},p_{+}^{n}):\\
&\left(
             \begin{array}{c}
              q_{-}^{n} \\
              p_{-}^{n} \\
             \end{array}
           \right)=\left(
             \begin{array}{c}
              q_{+}^{n} \\
              p_{+}^{n}+h\Upsilon(h \Omega)g(q_{-}^{n}) \\
             \end{array}
           \right),\\
&\textmd{3. }(q^{n+1},p^{n+1}) = \Phi_{h/2,\textmd{L}}
(q_{-}^{n},p_{-}^{n}): \\
&\left(
             \begin{array}{c}
              q^{n+1} \\
              p^{n+1} \\
             \end{array}
           \right)=\left(
                                                           \begin{array}{cc}
                                                            \cos(\frac{h\Omega}{2})& \Omega^{-1}\sin(\frac{h\Omega}{2}) \\
                                                             -\Omega \sin(\frac{h\Omega}{2}) & \cos(\frac{h\Omega}{2}) \\
                                                           \end{array}
                                                         \right)\left(
             \begin{array}{c}
               q_{-}^{n} \\
               p_{-}^{n}\\
             \end{array}
           \right).
\end{array}
\end{equation*}
It can be straightly verified that  this Strang splitting method is
identical to the ERKN integrator $\Phi_h$, i.e.,
 \begin{equation}
\Phi_h=\Phi_{h/2,\textmd{L}} \circ \Phi_{h,\textmd{NL}} \circ
\Phi_{h/2,\textmd{L}}
  \label{connetciton1}%
\end{equation}
if and only if the following conditions hold
 \begin{equation}
\frac{1}{2}\textmd{sinc}(\frac{1}{2}h \Omega)\Upsilon(h
\Omega)=\bar{b}_1(h \Omega),\ \ \cos(\frac{1}{2}h \Omega)\Upsilon(h
\Omega)=b_1(h \Omega).
  \label{connetciton10}%
\end{equation}
For general ERKN integrators, these two conditions may not  hold
simultaneously. However, it is noted that under the symmetry
condition \eqref{sym cond},   \eqref{connetciton10} is true. In this
particular case, the function $\Upsilon$ has the form
\begin{equation}  \label{Upsilon}
\Upsilon(h \Omega)=b_1(h \Omega)\cos^{-1}(\frac{1}{2}h
\Omega)=2\bar{b}_1(h \Omega)\textmd{sinc}^{-1}(\frac{1}{2}h
\Omega).\end{equation}

 On the other hand, if we consider the Strang splitting in
another way:
 \begin{equation*}
\hat{\Phi}_h=\Phi_{h/2,\textmd{NL}} \circ \Phi_{h,\textmd{L}} \circ
\Phi_{h/2,\textmd{NL}},
\end{equation*}
which yields a class of   trigonometric integrators
\begin{equation}\label{TI}
\begin{array}[c]{ll} &q^{n+1}=
\cos(h\Omega)q^{n}+h\textmd{sinc}(h\Omega)p^{n}+\frac{1}{2}h^2
 \textmd{sinc}(h\Omega)\Upsilon(h
\Omega)g(q^n),\\
 &p^{n+1}=-\Omega
 \sin(h\Omega)q^{n}+\cos(h\Omega)p^{n}+\frac{1}{2}h\big(
 \cos(h\Omega)\Upsilon(h \Omega)g(q^n)+\Upsilon(h
\Omega)g(q^{n+1})\big).
\end{array}
\end{equation}
It is noted that this is exactly a trigonometric integrator (a form
of (XIII.2.7)-(XIII.2.8) given on p.481 of \cite{hairer2006}) with
the following coefficients
\begin{equation}\label{phi deng}\begin{array}[c]{ll}\Phi(h \Omega)=1,\qquad \Psi(h \Omega)=
\textmd{sinc}(h\Omega)\Upsilon(h \Omega),\\
 \Psi_0(h \Omega)=
\cos(h\Omega)\Upsilon(h \Omega),\ \ \Psi_1(h \Omega)=\Upsilon(h
\Omega).\end{array}\end{equation} In terms of the symmetry condition
proposed in \cite{hairer2006}, it is clear that this trigonometric
integrator is symmetric.

On the basis of the above analysis, the following important
connection between symmetric ERKN integrators and symmetric
trigonometric integrators is obtained.
\begin{prop}\label{connection thm}
 For the one-stage explicit symmetric ERKN integrator \eqref{methods} and the symmetric trigonometric
integrator \eqref{TI}, they are related by
 \begin{equation}
\Phi_h=\Phi_{-h/2,\textmd{L}}\circ\Phi_{-h/2,\textmd{NL}}\circ\hat{\Phi}_h\circ
\Phi_{h/2,\textmd{NL}}\circ \Phi_{h/2,\textmd{L}}.
  \label{connetciton2}
\end{equation}
The function $\Upsilon$ appearing in the averaged equation and in
\eqref{phi deng} is determined by \eqref{Upsilon}. Moreover, it is
true that
 \begin{equation}\begin{aligned}
\underbrace{\Phi_h\circ \cdots \circ\Phi_h}_{n\
\textmd{times}}=&\Phi_{-h/2,\textmd{L}}\circ\Phi_{-h/2,\textmd{NL}}\circ\big(\underbrace{\hat{\Phi}_h\circ
\cdots \circ \hat{\Phi}_h}_{n\ \textmd{times}}\big)   \circ
\Phi_{h/2,\textmd{NL}}\circ \Phi_{h/2,\textmd{L}}\\
=&\Phi_{h/2,\textmd{L}}\circ\Phi_{h/2,\textmd{NL}}\circ\big(\underbrace{\hat{\Phi}_h\circ
\cdots \circ \hat{\Phi}_h}_{n-1\ \textmd{times}}\big)   \circ
\Phi_{h/2,\textmd{NL}}\circ \Phi_{h/2,\textmd{L}}.
  \label{connetciton3}%
\end{aligned}\end{equation}
\end{prop}

\subsection{Proof of Theorem \ref{thm: symetric Long-time
thm}} Since the long-time behaviour of symmetric  trigonometric
integrators  is well understood (see
\cite{Cohen15,Hairer00,hairer2006}), the   long-time near
conservation of total and oscillatory energy of one-stage explicit
symmetric ERKN integrators   can be derived  by using the relation
\eqref{connetciton3}. To be very precise, according to  the analysis
of \cite{Cohen15,Hairer00,hairer2006}, one has to verify the
following four key points.
\begin{itemize}
\item I. If $(q^0 ,p^0 )$  satisfies \eqref{energy condition}, the following initial value  $$(\tilde{q}^0,\tilde{p}^0 ):=\Phi_{h/2,\textmd{NL}}\circ \Phi_{h/2,\textmd{L}}(q^0,p^0 )$$ for $\big(\underbrace{\hat{\Phi}_h\circ
\cdots \circ \hat{\Phi}_h}_{n-1\ \textmd{times}}\big)$ satisfies the
finite-energy condition
\begin{equation}
\frac{1}{2} \norm{\tilde{p}^0}^2+\frac{1}{2} \norm{\Omega \tilde{q}^0}^2\leq E. \label{new IV energy condition}%
\end{equation}

\item II. For any $j\in \mathbb{N}^{+}$, $\big(\underbrace{\hat{\Phi}_h\circ
\cdots \circ \hat{\Phi}_h}_{j\ \textmd{times}}\big)   \circ
\Phi_{h/2,\textmd{NL}}\circ \Phi_{h/2,\textmd{L}}(q^0,p^0 )$ stays
in a compact set if $\big(\underbrace{\Phi_h\circ \cdots \circ
\Phi_h}_{j\ \textmd{times}}\big) (q^0,p^0 )$ does.


\item III. The two additional steps with $\Phi_{h/2,\textmd{L}}$ and $\Phi_{h/2,\textmd{NL}}$  only introduce an $\mathcal{O}(h) $ deviation in the
total and oscillatory   energy, provided that the corresponding
initial values $(\tilde{q},\tilde{p})$ are bounded and also $\Omega
\tilde{q}$ is bounded.

\item IV. The condition of \cite{Cohen15} on $\Phi$ and $\Psi$  has to be satisfied
for the choice \eqref{phi deng}.
\end{itemize}
In what follows, we show that the above four conditions are
completely true.
\begin{itemize}
\item For the point I, according to the splitting method, it is easy to
obtain that
\begin{equation*}
\begin{aligned} &\tilde{q}^0=\cos(\frac{h\Omega}{2})q^0+\Omega^{-1}\sin(\frac{h\Omega}{2})p^0,\\
&\tilde{p}^0=-\Omega \sin(\frac{h\Omega}{2})
q^0+\cos(\frac{h\Omega}{2})p^0+\frac{1}{2}h\Upsilon(\frac{h\Omega}{2})g(\tilde{q}^0)
.
\end{aligned}
\end{equation*} Thus \eqref{new IV energy condition} is clear  by
considering   \eqref{energy condition}.

\item  It follows from \eqref{connetciton2} that  \begin{equation*}\begin{aligned}
\big(\underbrace{\hat{\Phi}_h\circ \cdots \circ \hat{\Phi}_h}_{j\
\textmd{times}}\big)   \circ \Phi_{h/2,\textmd{NL}}\circ
\Phi_{h/2,\textmd{L}}(q^0,p^0 )\\
= \Phi_{h/2,\textmd{NL}}\circ\Phi_{h/2,\textmd{L}}\circ
\big(\underbrace{\Phi_h\circ \cdots \circ \Phi_h}_{j\
\textmd{times}}\big) (q^0,p^0 ),
\end{aligned}\end{equation*}
which gives the statement of II.

\item  The third result is clear in the light of the definitions
of  $\Phi_{h/2,\textmd{L}}$ and $\Phi_{h/2,\textmd{NL}}$.

\item  We remark that the long-term analysis of \cite{Hairer00,hairer2006}
requires $$|\Phi(h\omega)| \geq C
|\textmd{sinc}(\frac{1}{2}h\omega)|.$$ However, for the
trigonometric integrator \eqref{TI} with $\Phi=1$, it does   not
satisfy the above requirement. Therefore,  the long-term analysis of
\cite{Hairer00,hairer2006} can not be used here. Recently, the
authors in \cite{Cohen15} improved the analysis and presented a
long-term analysis of numerical integrators for oscillatory
Hamiltonian systems under minimal non-resonance conditions. A more
relaxed restriction on $\Phi$ was given there. In terms of that
restriction and the relations  \eqref{connetciton10} and
\eqref{connetciton2}, the statement of IV holds provided
\eqref{sigma bound} is true.
\end{itemize}

Based on the above analysis and the long-time result of symmetric
trigonometric integrators  given in
\cite{Cohen15,Hairer00,hairer2006}, Theorem   \ref{thm: symetric
Long-time thm} is immediately obtained.


\section{Proof for symplectic integrators} \label{sec: symplectic proof}
In this section, we are devoted to  the analysis of long-time
conservation of the total and oscillatory energy along symplectic
ERKN integrators.  To do this, we begin with defining the operators
\begin{equation*}
 \begin{array}[c]{ll}
\mathcal{L}_1(hD)=&
\big(e^{hD}-\cos(h\Omega)-h^2\Omega^2\sinc(h\Omega)\bar{b}_1(h\Omega)b^{-1}_1(h\Omega)\big)\\
&\big(\bar{b}_1(h\Omega)b^{-1}_1(h\Omega)(e^{hD}-\cos(h\Omega))+\sinc(h\Omega)\big)^{-1},\\
\mathcal{L}_2(hD)= &
\cos(c_1h\Omega)+c_1\sinc(c_1h\Omega)\mathcal{L}_1(hD),\\
\mathcal{L}(hD)=
&\big(e^{hD}-\cos(h\Omega)-\sinc(h\Omega)\mathcal{L}_1(hD)\big)
(\bar{b}_1(h\Omega)\mathcal{L}_2(hD))^{-1}
,\end{array}\end{equation*}
 where $D$ is the differential operator (see \cite{Hairer16}). The following property will be used in
this section.
\begin{prop}\label{lhd pro}
Under    the symplecticity condition \eqref{symple cond} with
$d_1=1$ of one-stage explicit ERKN integrators,  the operator
$\mathcal{L}(hD)$ can be simplified as
$$ \mathcal{L}(hD)=-2 h \Omega(\cos(h\Omega) - \cosh(hD))
\csc(h\Omega).$$ Moreover, the Taylor expansions of the operator
$\mathcal{L}(hD)$ are
\begin{equation*}
\begin{aligned}
\mathcal{L}(hD)=&  2h \Omega  \tan(\frac{1}{2}h\Omega)+h \Omega
\csc(h\Omega)
(hD)^2+\cdots,\\
\mathcal{L}(hD+\mathrm{i}h\omega)=& \left(
                                              \begin{array}{cc}
                                                -4 \sin^2(\frac{1}{2}h\omega)) &  \\
                                                 & 0 \\
                                              \end{array}
                                            \right)
+\left(
                                              \begin{array}{cc}
                                                2\sin(h\omega) &  \\
                                                 & 2 h\omega\\
                                              \end{array}
                                            \right)(\mathrm{i}hD)+\cdots,\\
\mathcal{L}(hD+\mathrm{i}kh\omega)=& -4h \Omega
\csc(h\Omega)\sin\big(\frac{kh\omega
I+h\Omega}{2}\big)\sin\big(\frac{kh\omega
I-h\Omega}{2}\big)\\
&+2h \Omega \csc(h\Omega)\sin(kh\omega)(\mathrm{i}hD)+\cdots, \
\textmd{where}\ \abs{k}>1.
\end{aligned}
\end{equation*}
\end{prop}

\subsection{Modulated Fourier expansion} \label{sec:Analysis of the
methods}
 We now  derive the modulated Fourier expansions of symplectic ERKN integrators
and present the bounds of the modulated Fourier functions.
\begin{mytheo}\label{energy thm}
Under the conditions of Theorem \ref{thm: symplec Long-time thm} and
for $0 \leq t=nh \leq T$, the numerical solution of the one-stage
explicit symplectic ERKN integrator \eqref{methods} admits the
following modulated Fourier expansions
\begin{equation}
\begin{aligned} &q^{n}= \sum\limits_{|k|<N} \mathrm{e}^{\mathrm{i}k\omega t}\zeta_h^k(t)+R_{h,N}(t),\\
&p^{n}= \sum\limits_{|k|<N} \mathrm{e}^{\mathrm{i}k\omega t}\eta_h^k(t)+S_{h,N}(t),\\
\end{aligned}
\label{MFE-ERKN}%
\end{equation}
where the remainder terms are bounded by
\begin{equation}
 R_{h,N}(t)=\mathcal{O}(th^{N}),\ \ \ \  S_{h,N}(t)=\mathcal{O}(th^{N-1}).\\
\label{remainder}%
\end{equation}
The coefficient functions $\zeta_h^k$  as well as all their
derivatives are bounded by
\begin{equation}
\begin{array}{ll}\ddot{\zeta}^0_{h,1}=\mathcal{O}(1),\quad  \
\zeta_{h,1}^1=\mathcal{O}(h^{2}), \ \
\zeta_{h,1}^k=\mathcal{O}(h^{\abs{k}+1}), \\
\zeta^0_{h,2}=\mathcal{O}(h^{\frac{3}{2}}),\ \
\dot{\zeta}_{h,2}^1=(h), \quad \ \ \
\zeta_{h,2}^k=\mathcal{O}(h^{\abs{k}+1}),
\end{array}\label{coefficient func1-num}
\end{equation}
for $\abs{k}>1$ and further bounded by
\begin{equation}\label{coefficient func1-1-num}%
 \zeta^0_{h,1}=\mathcal{O}(1), \quad  \ \
\zeta_{h,2}^1=\mathcal{O}(h).
\end{equation}
Moreover, we have the following results for coefficient functions
$\eta_h^k$
\begin{equation}
\begin{array}{ll}
\eta_{h,1}^0=\dot{\zeta}_{h,1}^0+\mathcal{O}(h),\  \
\eta_{h,1}^1=\mathcal{O}(h),\qquad \qquad \
\eta_{h,1}^k=\mathcal{O}(h^{\abs{k}}),
\\   \eta_{h,2}^0=\mathcal{O}(h^{\frac{1}{2}}), \qquad \quad \eta_{h,2}^1=
\mathrm{i}\omega \zeta_{h,2}^1+\mathcal{O}(h), \ \
\eta_{h,2}^k=\mathcal{O}(h^{\abs{k}}),
\end{array}
\label{rea20-num}%
\end{equation}
where  $\abs{k}>1$.
 Moreover, we have
$\zeta^{-k}=\overline{\zeta^{k}}$ and
$\eta^{-k}=\overline{\eta^{k}}$. The constants symbolized by the
notation are independent of $h$ and $\omega$, but depend on the
constants  from Assumption \ref{ass} and the final time $T$.
\end{mytheo}
\textbf{Proof.}  The proof follows   that used in the modulated
Fourier expansions  of previous publications (see
\cite{Cohen15,Cohen05,Hairer00,hairer2006}) but with some novel
adaptations for non-symmetric methods. The proof of this theorem
does not rely on the symmetry of the methods, which is a main
conceptual difference in comparison with that in
\cite{Cohen15,Cohen05,Hairer00,hairer2006}.


We will prove that there exist  two functions
\begin{equation}
\begin{aligned} &q_{h}(t)= \sum\limits_{|k|<N} \mathrm{e}^{\mathrm{i}k\omega
t}\zeta_h^k(t),\
 \ p_{h}(t)= \sum\limits_{|k|<N} \mathrm{e}^{\mathrm{i}k\omega t}\eta_h^k(t)\\
\end{aligned}
\label{MFE-1}%
\end{equation}
with smooth (in the sense that all their derivatives are bounded
independently of $h$ and $\omega$) coefficients $ \zeta_h^k,
\eta_h^k$,  such that, for $t=nh$,
\begin{equation*}
\begin{aligned} &q^{n}=q_{h}(t)+\mathcal{O}(h^{N}),\qquad p^{n}= p_{h}(t)+\mathcal{O}(h^{N-1}).\\
\end{aligned}
\end{equation*}

\vskip1mm \textbf{Construction of the coefficients functions.}
According to the second and third formulae of \eqref{methods}, it is
arrived that
\begin{equation*}
\begin{aligned} &q^{n+1}-\big(\cos(h\Omega)+h^2\Omega^2\sinc(h\Omega)\bar{b}_1(h\Omega)b^{-1}_1(h\Omega)\big)q^{n}\\
=&h\bar{b}_1(h\Omega)b^{-1}_1(h\Omega)p^{n+1}+h\big(\sinc(h\Omega)-\bar{b}_1(h\Omega)b^{-1}_1(h\Omega)\cos(h\Omega)\big)p^{n}.
\end{aligned}
\end{equation*}
Comparing the coefficients of $\mathrm{e}^{\mathrm{i}k\omega t}$ and
considering the definition of $\mathcal{L}_1$ implies the
relationship between $\zeta_h^k$ and $\eta_h^k$ as follows:
\begin{equation}\label{MFE-zetaeta}%
\begin{aligned}& \mathcal{L}_1(hD)\zeta_h^0=h\eta_h^0,\\
&\mathcal{L}_1(hD+\mathrm{i}kh\omega)\zeta_h^k=h\zeta_h^k\qquad
\textmd{for}\ \ \  \abs{k}>0.
\end{aligned} %
\end{equation}

  For the first formula of the ERKN integrator \eqref{methods}, we look for the  function
\begin{equation}
\begin{aligned} & \tilde{q}_h(t):=\sum\limits_{|k|<N}
\mathrm{e}^{\mathrm{i}k\omega t}\xi_h^k(t)
\end{aligned}
\label{MFE-2}%
\end{equation}
as the modulate Fourier expansion of   $Q^{n+\frac{1}{2}}$ at
$t=(n+c_1)h$. Inserting \eqref{MFE-1}-\eqref{MFE-2} into the first
formula of  \eqref{methods} and comparing the coefficients of
$\mathrm{e}^{\mathrm{i}k\omega t}$ yields
\begin{equation}\label{MFE-3}%
\begin{aligned}& \mathcal{L}_2(hD)\zeta_h^0=\xi_h^0,\\
&\mathcal{L}_2(hD+\mathrm{i}kh\omega)\zeta_h^k=\xi_h^k \qquad
\textmd{for}\ \ \ \abs{k}>0,
\end{aligned} %
\end{equation}
which gives the relationship between the coefficient functions
$\xi_h^k$ and $\zeta_h^k$.

We insert the above expansions into the second equation of
\eqref{methods}, expand  the nonlinear function into its Taylor
series and compare the coefficients of
$\mathrm{e}^{\mathrm{i}k\omega t}$. We then obtain
\begin{equation*}
\begin{aligned}& \mathcal{L}_2(hD)\zeta_h^0=h^2\Big(g(\zeta_h^0)+
\sum\limits_{s(\alpha)=0}\frac{1}{m!}g^{(m)}(\zeta_h^0)(
\zeta_h)^{\alpha}\Big),\\
& \mathcal{L}_2(hD+\mathrm{i}kh\omega)\zeta_h^k=h^2
\sum\limits_{s(\alpha)=k}\frac{1}{m!}g^{(m)}(\zeta_h^0)(
\zeta_h)^{\alpha},
\end{aligned} %
\end{equation*}
where $\abs{k}>0$, the sum ranges over $m\geq0$, the multi-indices
$\alpha=(\alpha_1,\ldots,\alpha_m)$ with integer $\alpha_i$
satisfying $0<|\alpha_i|<N$   have  a given sum
$s(\alpha)=\sum\limits_{j=1}^{m}\alpha_j,$ and $(\zeta_h)^{\alpha}$
is an abbreviation for the $m$-tuple
$(\zeta^{\alpha_1}_h,\ldots,\zeta^{\alpha_m}_h)$.

Comparing the dominate terms in the relations for the coefficients
functions $\zeta^k_h$ motivates the following ansatz of the
modulated Fourier functions:
\begin{equation}\label{ansatz}%
\begin{array}{ll}
\ddot{\zeta}^0_{h,1}=\big(\mathcal{G}_{10}(\cdot)+\cdots\big),\
&\zeta^0_{h,2}= \frac{ h^2 \cos(\frac{1}{2}h\omega)}{ 2h \omega  \sin(\frac{1}{2}h\omega)}\mathcal{G}_{20}(\cdot)+\cdots,\\
\zeta_{h,1}^1=\frac{ h^2}{-4 \sin^2(\frac{1}{2}h\omega))}
\big(\mathcal{F}^1_{10}(\cdot)+\cdots\big),\
&\dot{\zeta}_{h,2}^1=\frac{
 h^2}{2\mathrm{i} h^2  \omega }
\big(\mathcal{F}^1_{20}(\cdot)+\cdots\big),
\\
\zeta^{k}=\frac{ h^2\sin(h\Omega)}{-4h \Omega
\sin\big(\frac{kh\omega I+h\Omega}{2}\big)\sin\big(\frac{kh\omega
I-h\Omega}{2}\big)} \big(\mathcal{F}^{k}_{0}(\cdot)+\cdots\big)
&\textmd{for}\ \abs{k}>1,
\end{array} %
\end{equation}
where  the dots  stand  for power series in $\sqrt{h}$. Since the
series in the ansatz usually diverge, in this paper we
 truncate them after the $\mathcal{O}(h^{N+1})$ terms (see \cite{Hairer00}).

\vskip1mm \textbf{Initial values.} In this part, we derive the
initial values
 for the differential equations appearing in the ansatz
\eqref{ansatz}.  On the basis of  the conditions $p_{h}(0)=p^0$ and
 $q_{h}(0)=q^0$,
it can be deduced that
\begin{equation}\label{Initial values-1}%
\begin{aligned}
&p^0_{1}=\eta^0_{h,1}(0)+\mathcal{O}(h)=\dot{\zeta}^0_{h,1}(0)+\mathcal{O}(h),\ &p^0_{2}=2\mathrm{Re}(\eta^1_{h,2}(0))+\mathcal{O}(h^{\frac{1}{2}}),\\
&q^0_{1}=\zeta^0_{h,1}(0)+\mathcal{O}(h^2),\
&q^0_{2}=2\mathrm{Re}(\zeta^1_{h,2}(0))+\mathcal{O}(h^{\frac{3}{2}}).\end{aligned} %
\end{equation}
On the other hand, we have
\begin{equation*}
\begin{aligned}&p_{h,1}(h)=p^1_{1},\ \  p_{h,2}(h)=p^1_{2},\ \
q_{h,1}(h)=q^1_{1},\ \ q_{h,2}(h)=q^1_{2}.\end{aligned}
\end{equation*}
 It follows from the scheme of the method \eqref{methods} that
\begin{equation}\label{Initial values-f1}
q_2^{1}-\cos(h\omega)q_2^{0}=h
\mathrm{sinc}(h\omega)p_2^{0}+h^2\bar{b}_1(h\omega)g_2(Q^{n+c_1}),\end{equation}
where we have used the notation:
$$
 g(Q^{n+c_1})=\big(g_1(Q^{n+c_1})^{\intercal},g_2(Q^{n+c_1})^{\intercal}\big)^{\intercal}.$$
By computing $q_2^{1}-\cos(h\omega)q_2^{0}$, we get
\begin{equation*}
\begin{aligned}
&q_2^{1}-\cos(h\omega)q_2^{0}=q_{h,2}(h)-\cos(h\omega)q_{h,2}(0)\\
=&\sum\limits_{|k|<N} \big(\mathrm{e}^{\mathrm{i}k\omega
h}\zeta_{h,2}^k(h)-\cos(h\omega)
\zeta_{h,2}^k(0)\big)\\
=&\zeta^0_{h,2}(h)+ \mathrm{e}^{\mathrm{i}\omega h}\zeta_{h,2}^1(h)+
\mathrm{e}^{-\mathrm{i}\omega
h}\zeta_{h,2}^{-1}(h)\\
&-\cos(h\omega)\Big(\zeta^0_{h,2}(0)+
\zeta_{h,2}^1(0)+\zeta_{h,2}^{-1}(0)\Big)+\mathcal{O}(h^2).
\end{aligned}
\end{equation*}
Expanding the functions $\zeta^0_{h,2}(h),\ \zeta_{h,2}^{1}(h),\
\zeta_{h,2}^{-1}(h)$ at $h=0$ yields
\begin{equation*}
\begin{aligned}
q_2^{1}-\cos(h\omega)q_2^{0}=&(1-\cos(h\omega))\zeta^0_{h,2}(0)+\mathrm{i}\sin(h\omega)
(\zeta_{h,2}^1(0)-\zeta_{h,2}^{-1}(0))+\mathcal{O}(h^2)\\
=&2\sin^2(h\omega/2)\zeta^0_{h,2}(0)+\mathrm{i}\sin(h\omega)
(\zeta_{h,2}^1(0)-\zeta_{h,2}^{-1}(0))+\mathcal{O}(h^2)\\
=&\mathrm{i}\sin(h\omega)
(\zeta_{h,2}^1(0)-\zeta_{h,2}^{-1}(0))+\mathcal{O}(h^2),
\end{aligned}
\end{equation*}
where we have used the result of $\zeta^0_{h,2}$ presented in
\eqref{ansatz}. Now  the formula \eqref{Initial values-f1} becomes
\begin{equation*}
\mathrm{i}\sin(h\omega) (\zeta_{h,2}^1(0)-\zeta_{h,2}^{-1}(0))  =h
\mathrm{sinc}(h\omega)p_2^{0}+h^2\bar{b}_1(h\omega)g_2(Q^{n+c_1})+\mathcal{O}(h^2),
\end{equation*}
which  confirms that
\begin{equation}\label{Initial values-2}2\mathrm{Im}(\zeta_{h,2}^1(0))=\omega^{-1}p_2^{0}+\mathcal{O}(\omega^{-1}).\end{equation}
Therefore, using the implicit function theorem, conditions
\eqref{Initial values-1} and \eqref{Initial values-2} yield  the
desired initial values $\zeta^0_{h,1}(0),\ \dot{\zeta}^0_{h,1}(0),\
\zeta_{h,2}^1(0)$ for the differential equations appearing in the
ansatz \eqref{ansatz}.

\vskip1mm

\textbf{Bounds of the coefficients functions.} Based on the ansatz
and Assumption \ref{ass}, it is easy to get
  the bounds \eqref{coefficient func1-num}.
 From the condition
\eqref{energy condition}, it follows that $q^0_2=\mathcal{O}(h).$
Then by this result, \eqref{Initial values-2} and the fourth formula
of \eqref{Initial values-1}, we get
$\zeta_{h,2}^1(0)=\mathcal{O}(h).$ This implies that $
\zeta^1_{h,2}(t)=\mathcal{O}(h)$ by considering
$\dot{\zeta}^1_{h,2}=\mathcal{O}(\omega^{-1}).$ Similarly, one
arrives at that $\zeta^0_{h,1}(0)=\mathcal{O}(1),
\dot{\zeta}^0_{h,1}(0)=\mathcal{O}(1)$ and then $
\zeta^0_{h,1}(t)=\mathcal{O}(1).$ Therefore \eqref{coefficient
func1-1-num} is obtained. The bounds \eqref{rea20-num} can be
derived easily by considering \eqref{MFE-zetaeta} and the bounds of
$\zeta_h^k$.

\vskip1mm

 \textbf{Defect.} As the last part of the proof, we analyze the defect. Firstly, define the components of the defect, for
 $t=nh$,
\begin{equation*}
\begin{aligned}
 d_1(t+h)=&q_{h,1}(t+h)-q_{h,1}(t)-hp_{h,1}(t)-h^2\bar{b}_1(0)g_1(\tilde{q}_{h}(t+c_1 h)),\\
 d_2(t+h)= &p_{h,1}(t+h)-p_{h,1}(t)-hb_1(0)g_1(\tilde{q}_{h}(t+c_1 h)),\\
 d_3(t+h)=&q_{h,2}(t+h) -\cos(h\omega)q_{h,2}(t)-h \mathrm{sinc}(h\omega)p_{h,2}(t)\\
 &-h^2\bar{b}_1(h\omega)g_2(\tilde{q}_{h}(t+c_1 h)),\\
 d_4(t+h)=&p_{h,2}(t+h)-\cos(h\omega)p_{h,2}(t)+\omega \sin(h\omega)q_{h,2}(t)\\
 &-hb_1(h\omega)g_2(\tilde{q}_{h}(t+c_1 h)).
\end{aligned}
\end{equation*}
By  the definition of  the coefficient functions $\zeta^k_{h},\
\eta^k_{h},\ \xi_h^k$, the following  results are  true
$$d_1=\mathcal{O}(h^{N+1}),\ d_2=\mathcal{O}(h^{N+1}),\
d_3=\mathcal{O}(h^{N+1}),\ d_4=\mathcal{O}(h^{N+1}).$$

We are now in a position to estimate the remainders
\eqref{remainder}. To do this, we begin with defining
$$R_{k}^n=p_k^n-p_{h,k}(t),\quad S_k^n=q_k^n-q_{h,k}(t)$$  for $k=1,2,$  and the norm
$$\norm{(S_1^n,R_1^n,S_2^n,R_2^n)}_{\ast}=\norm{(S_1^n,R_1^n,\omega S_2^n,R_2^n)}.$$
We first estimate the remainder for the difference
$Q^{n+c_1}-\tilde{q}_{h}(t+c_1 h).$ By   the first two equations in
\eqref{methods} and the first equation in \eqref{MFE-2}, if the
function $g$ satisfies a Lipschitz condition, one obtains
$$\norm{Q^{n+c_1}-\tilde{q}_{h}(t+c_1 h)}\leq\norm{S^n}+c_1h\norm{R^n}\triangleq \alpha^n.$$
 For the remainders,  on noticing  the scheme of the ERKN integrators, we then have
\begin{equation}\label{relaSR}
\begin{aligned}
\norm{(S_1^{n+1},R_1^{n+1},S_2^{n+1},R_2^{n+1})}_{\ast}\leq&
\norm{(S_1^n,R_1^n,S_2^n,R_2^n)}_{\ast}+h\kappa_1
\norm{(S_1^n,R_1^n,S_2^n,R_2^n)}_{\ast}\\
&+h\kappa_2 \alpha^n+\kappa_3 h^{N+1},
\end{aligned}
\end{equation}
where $\kappa_1,\ \kappa_2,\ \kappa_3$ are constants. From the
definition of the initial values,  it follows that $\norm{
(S_1^0,R_1^0,S_2^0,R_2^0) }_{\ast}=\mathcal{O}(h^{N+1}).$ By using
the relation \eqref{relaSR} repeatedly, we obtain the following
estimate for the remainders
$\norm{(S_1^n,R_1^n,S_2^n,R_2^n)}_{\ast}\leq Cnh^{N+1}, $ which
yields \eqref{remainder}.
\\

We complete the proof of this theorem.

\subsection{The first almost-invariant}
 Let
\begin{equation*}
\begin{aligned}
&\vec{\zeta}=\big(\zeta^{-N+1}_h,\cdots,\zeta^{-1}_h,
\zeta^{0}_h,\zeta^{1}_h,\cdots,\zeta^{N-1}_h\big),\\
&\vec{\eta}=\big(\eta^{-N+1}_h,\cdots,\eta^{-1}_h,
\eta^{0}_h,\eta^{1}_h,\cdots,\eta^{N-1}_h\big).
\end{aligned}
\end{equation*} We
have the  first almost-invariant of the modulated Fourier functions
as follows.
\begin{mytheo}\label{first invariant thm}
Under the conditions of Theorem \ref{energy thm}, there exists a
function $\widehat{\mathcal{H}}[\vec{\zeta},\vec{\eta}]$ such that
the coefficient functions of   the modulated Fourier expansion of
symplectic ERKN integrators  satisfy
\begin{equation}
\widehat{\mathcal{H}}[\vec{\zeta},\vec{\eta}](t)=\widehat{\mathcal{H}}[\vec{\zeta},\vec{\eta}](0)+\mathcal{O}(th^{N})
\label{HH}%
\end{equation}
for $0\leq t\leq T.$ Moreover,   this  can  be expressed as
\begin{equation}\begin{aligned}
\widehat{\mathcal{H}}[\vec{\zeta},\vec{\eta}]=\frac{1}{2}
(\eta^0_{h,1})^\intercal \eta^0_{h,1}
+2\omega^2\big(\zeta_{h,2}^{-1}\big)^\intercal\zeta_{h,2}^{1}
+U(\zeta^0_h(t))+\mathcal{O}(h^2).
\label{HH def}%
\end{aligned}
\end{equation}
\end{mytheo}
\textbf{Proof.} With Theorem \ref{energy thm}, we obtain
\begin{equation*}
\begin{aligned}
& \mathcal{L}(hD) q_{h}(t)=h^2 g( q_{h}(t))+\mathcal{O}(h^{N+2}),
\end{aligned}
\end{equation*}
where   the following denotations are used:
\begin{equation*}
\begin{aligned}q_{h}(t)=\sum\limits_{ |k|<N}q^k_{h}(t) \  \   \
\textmd{with}\ \ \ q^k_{h}(t)=\mathrm{e}^{\mathrm{i}k\omega
t}\zeta_h^k(t).
\end{aligned}
\end{equation*}
 Considering the  definition of $ q_h$ and
comparing the coefficients of $\mathrm{e}^{\mathrm{i}k\omega t}$
yields the resulting equations in terms of $ q_h^k:$
\begin{equation}
\begin{aligned}
& \mathcal{L}(hD) q^k_{h}(t)=-h^2
\nabla_{q^{-k}}\mathcal{U}(\vec{q}_{h}(t))+\mathcal{O}(h^{N+2}),
\end{aligned}
\label{methods-inva-nnew}%
\end{equation}
where $\mathcal{U}(\vec{q}_{h}(t))$ is defined as
\begin{equation}
\begin{aligned}
&\mathcal{U}(\vec{q}_{h}(t))=U(q^0_{h}(t))+
\sum\limits_{s(\alpha)=0}\frac{1}{m!}U^{(m)}(q^0_{h}(t))
(q_{h}(t))^{\alpha},
\end{aligned}
\label{newuu}%
\end{equation}
and $\vec{q}_{h}(t)$ is given by
\begin{equation*}
\begin{aligned}
\vec{q}_{h}(t)=\big(q^{-N+1}_{h}(t),\ldots,
q^{-1}_{h}(t),q^{0}_{h}(t),q^{1}_{h}(t),\ldots,q^{N}_{h}(t)\big).
\end{aligned}
\end{equation*}
Multiplying the equation \eqref{methods-inva-nnew} with $
\big(\dot{q}_h^{-k}\big)^\intercal$ and summing up conforms
\begin{equation*}
\begin{aligned} & \frac{1}{h^2}\sum\limits_{|k|<N}
 \big(\dot{q}_h^{-k}\big)^\intercal  \mathcal{L}(hD)
q^k_{h}+\frac{d}{dt}\mathcal{U}(\vec{q}_{h})=\mathcal{O}(h^{N}).
\end{aligned}
\end{equation*}
 We switch to the
quantities $\zeta_h^k(t)$ and get the equivalent relation
\begin{equation}
\begin{aligned}
 \mathcal{O}(h^{N}) =&  \frac{1}{h^2}\sum\limits_{|k|<N}
 \big(\dot{\zeta}_h^{-k}-\textmd{i} k \omega \zeta_h^{-k}\big)^\intercal
\mathcal{L}(hD+\mathrm{i}hk\omega) \zeta_h^k+\frac{d}{dt}\mathcal{U}(\vec{\zeta})\\
= &\frac{1}{h^2}\sum\limits_{|k|<N}
 \big(\dot{\bar{\zeta}}_h^{k}-\textmd{i} k \omega \bar{\zeta}_h^{k}\big)^\intercal
\mathcal{L}(hD+\mathrm{i}hk\omega)
\zeta_h^k+\frac{d}{dt}\mathcal{U}(\vec{\zeta}).
\end{aligned}
\label{duu-new1}%
\end{equation}

 With the Taylor expansions of $\mathcal{L}(hD)$ and
 $\mathcal{L}(hD+\textmd{i}kh\omega)$ given in Proposition \ref{lhd pro}
 and   the ``magic formulas" on p. 508 of \cite{hairer2006}, we
 know that the following   part  appearing in \eqref{duu-new1}  is a total derivative
\begin{equation*}
\begin{aligned}
\sum\limits_{|k|<N}\big(\dot{\bar{\zeta}}_h^{k}-\textmd{i} k \omega
\bar{\zeta}_h^{k}\big)^\intercal \mathcal{L}(hD+\mathrm{i}hk\omega)
\zeta_h^k.
\end{aligned}
\end{equation*}
 Thus, the right-hand side  of \eqref{duu-new1} is a total derivative.
 Therefore, by \eqref{duu-new1} and the above analysis,  it can be confirmed that there exists a function $\widehat{\mathcal{H}}$ such
 that
$\frac{d}{dt}\widehat{\mathcal{H}}[\vec{\zeta},\vec{\eta}](t)=\mathcal{O}(h^{N})$
and an integration yields the statement \eqref{HH} of the theorem.

 By the previous analysis,  the bounds of Theorem \ref{energy thm}, Proposition \ref{lhd pro} and the analysis given in \cite{hairer2006}, we
obtain the construction of $\widehat{\mathcal{H}}$ as
\begin{equation*}
\begin{aligned}\widehat{\mathcal{H}}[\vec{\zeta},\vec{\eta}]=&\frac{1}{2}
(\dot{\zeta}^0_{h})
 ^\intercal \frac{h\Omega}{\sin(h\Omega)}\dot{\zeta}^0_{h}
+2\frac{\omega}{2h^2} 2h^2 \omega \big(\zeta_{h,2}^{-1}\big)^\intercal\zeta_{h,2}^{1}+U(\zeta^0_h)+\mathcal{O}(h^2)\\
=&\frac{1}{2} (\dot{\zeta}^0_{h,1})
 ^\intercal \dot{\zeta}^0_{h,1}
+ 2 \omega^2 \big(\zeta_{h,2}^{-1}\big)^\intercal\zeta_{h,2}^{1}+U(\zeta^0_h)+\mathcal{O}(h^2)\\
=&\frac{1}{2} (\eta^0_{h,1})
 ^\intercal \eta^0_{h,1}
+2 \omega^2
\big(\zeta_{h,2}^{-1}\big)^\intercal\zeta_{h,2}^{1}+U(\zeta^0_h)+\mathcal{O}(h^2),
\end{aligned}
\end{equation*}
where the fact that
$\dot{\zeta}^0_{h,1}=\eta^0_{h,1}+\mathcal{O}(h)$
 from  \eqref{rea20-num} is used. The proof is complete.

\subsection{The second  almost-invariant}
\begin{mytheo}\label{second invariant thm}
Under the conditions of Theorem \ref{first invariant thm}, there
exists a function $\widehat{\mathcal{I}}[\vec{\zeta},\vec{\eta}]$
such that  the coefficient functions of   the modulated Fourier
expansion of symplectic ERKN integrators satisfy
\begin{equation}
\widehat{\mathcal{I}}[\vec{\zeta},\vec{\eta}](t)=\widehat{\mathcal{I}}[\vec{\zeta},\vec{\eta}](0)+\mathcal{O}(th^{N})
\label{II}%
\end{equation}
for $0\leq t\leq T.$ Moreover, this  can be expressed as
\begin{equation}\begin{aligned}
\widehat{\mathcal{I}}[\vec{\zeta},\vec{\eta}] =&2 \omega^2
\big(\zeta_{h,2}^{-1}\big)^\intercal\zeta_{h,2}^{1}+U(\zeta^0_h)+\mathcal{O}(h^2).
\label{II def}
\end{aligned}
\end{equation}
\end{mytheo}
\textbf{Proof.}
  Define  the vector function
$\vec{q}(\lambda,t)$ of $\lambda$ as
$$\vec{q}(\lambda,t)=\big( \mathrm{e}^{\mathrm{i}(-N+1)\lambda
\omega}q^{-N+1}_h(t),\cdots,
q^{0}_h(t),\cdots,\mathrm{e}^{\mathrm{i}(N-1)\lambda
\omega}q^{N-1}_h(t)\big).$$ Then it follows from the definition
\eqref{newuu}
  that $\mathcal{U}( \vec{q}(\lambda,t))$ does not depend on
$\lambda$. Thus, its derivative with respect to $\lambda$  yields
\begin{equation*}
\begin{aligned}0=&\frac{d}{d\lambda}\mathcal{U}( \vec{q}(\lambda,t))
=\sum\limits_{|k|<N}\mathrm{i}k\omega\mathrm{e}^{\mathrm{i}k\lambda
\omega} (q^{k}_h(t))^\intercal \nabla_{q^{k}}\mathcal{U}(
\vec{q}(\lambda,t)).\end{aligned}
\end{equation*}
 Letting $\lambda=0$ yields
$ \sum\limits_{|k|<N}\mathrm{i}k\omega (q^{k}_h(t))^\intercal
\nabla_{q^{k}}\mathcal{U}( \vec{q}_h(t))=0.$ Therefore, one obtains
\begin{equation*}
\begin{aligned}
0=&\sum\limits_{|k|<N}\mathrm{i}k\omega (q^{-k}_h(t))^\intercal
\nabla_{q^{-k}}\mathcal{U}( \vec{q}_h(t))\\
 =&\frac{\mathrm{i}
\omega}{h^2}\sum\limits_{|k|<N}k  \big(q^{-k}_{h}(t)\big)^\intercal
\mathcal{L}(hD) q^k_{h}(t)
 +\mathcal{O}(h^{N}).
\end{aligned}
\end{equation*}
Rewritten in  the   $\zeta_h^k(t)$ variables, this becomes
\begin{equation*}
\begin{aligned}
&\frac{\mathrm{i} \omega}{h^2}\sum\limits_{|k|<N}k  \big(
\zeta_h^{-k} (t)\big)^\intercal \mathcal{L}(hD+\textmd{i}k\omega h)
\zeta^k_{h}(t)  =\mathcal{O}(h^{N}).
\end{aligned}
\end{equation*}
As in the proof of Theorem \ref{first invariant thm}, the left-hand
expression of this equation can be written as the time derivative of
a function. Therefore,     we get
$\frac{d}{dt}\widehat{\mathcal{I}}[\vec{\zeta},\vec{\eta}](t)=\mathcal{O}(h^{N})$
and an integration yields statement \eqref{II} of the theorem.

According to the above analysis and the bounds of Theorem
\ref{energy thm}, the construction of
$\widehat{\mathcal{I}}[\vec{\zeta},\vec{\eta}](t)$  is obtained,
 which concludes the proof.

\subsection{Long-time near-conservation of total and oscillatory
energy}\label{sec:Long-time near-conservation} Theorem \ref{thm:
symplec Long-time thm} will be proved in this subsection. Before
that, we give the following theorem.
\begin{mytheo}\label{HHthm} Under the
conditions of Theorem \ref{first invariant thm}, it holds that
\begin{equation*}
\begin{aligned}
&\widehat{\mathcal{H}}[\vec{\zeta},\vec{\eta}](nh)=H(q_n,p_n)+\mathcal{O}(h),\\
&\widehat{\mathcal{I}}[\vec{\zeta},\vec{\eta}](nh)=I(q_n,p_n)+\mathcal{O}(h),
\end{aligned}
\end{equation*}
where the constants symbolized by $\mathcal{O}$  depend on $N,\ T$
and the constants in the assumptions,  and  $\sigma(h\omega)$ is
given by \eqref{sigma}.
\end{mytheo}
\textbf{Proof.} In terms of  the analysis given   in Section XIII of
\cite{hairer2006} and the bounds presented in Theorem \ref{energy
thm}, one arrives at
 \begin{equation}\label{HIE}\begin{aligned} &H(q^n,  p^n)=\frac{1}{2}
 (\eta^0_{h,1})
 ^\intercal\eta^0_{h,1} +2\omega^2
\big(\zeta_{h,2}^{-1}\big)^\intercal\zeta_{h,2}^{1}
+U( \zeta^0_h)+\mathcal{O}(h),\\
&I(q^n, p^n)=2\omega^2
\big(\zeta_{h,2}^{-1}\big)^\intercal\zeta_{h,2}^{1} +\mathcal{O}(h).
\end{aligned}\end{equation}
A comparison between \eqref{HH def}, \eqref{II def}  and \eqref{HIE}
gives the stated relations of this theorem. \hfill

 On the basis of previous analysis given in this section and
following the approach used in Chapter XIII of \cite{hairer2006},
Theorem \ref{thm: symplec Long-time thm} is easily proved by
patching together the local near-conservation result.

\section{Conclusions} \label{sec:conclusions}

In this paper, the long-time total and oscillatory energy
conservation behaviour of one-stage explicit ERKN integrators was
studied when applied to  highly oscillatory Hamiltonian systems. It
turned out that a good long-time energy conservation holds not only
for symmetric integrators but also for symplectic integrators. A
relationship between symmetric ERKN integrators and symmetric
trigonometric integrators was established and on the basis of which,
the long-time conservation for symmetric ERKN integrators was
proved. In order to show the result for symplectic ERKN integrators,
modulated Fourier expansions    were developed  with some novel
adaptations for symplectic (not necessarily symmetric) methods.
Using this technique,  we proved that one-stage symplectic ERKN
integrators  have two almost-invariants and approximately conserve
the  energy over long times.

This is a preliminary  research on the long-time behaviour of ERKN
integrators for highly oscillatory Hamiltonian systems  and the
authors are clearly aware that there are still  some issues which
will  be further considered.
\begin{itemize}\itemsep=-0.2mm

\item  The long-time behaviour of symmetric or symplectic ERKN
integrators for multi-frequency highly oscillatory Hamiltonian
systems will be discussed in another work.

\item    The     energy
conservation behaviour of symmetric or symplectic ERKN integrators
in other ODEs such as Hamiltonian systems with a solution-dependent
high frequency or without   any non-resonance condition  will also
be considered.

\item Another issue for future exploration is the near-conservation  of
energy, momentum and actions along symmetric or symplectic ERKN
integrators of Hamiltonian wave equations.

\item We only consider  one-stage  symmetric or symplectic ERKN integrators in this paper.
The extension of this paper's analysis to   higher-stage ERKN
integrators is not obvious since  there is the technical difficulty
which needs to be overcome. This issue will be considered in future
investigations.

\end{itemize}

\section*{Acknowledgements}

We are grateful to Professor Christian Lubich for his helpful
comments and discussions on the topic of modulated Fourier
expansions. We also thank  Ludwig Gauckler for drawing our attention
to the   connection between symmetric ERKN methods and symmetric
trigonometric integrators.


\end{document}